\documentclass[11pt]{article}
\usepackage{amsmath,amsfonts,amsthm,amscd,amssymb,graphicx, color}

\usepackage[utf8]{inputenc}
\usepackage{amsbsy}
\usepackage{graphicx}
\usepackage[text={7.0in,9.5in},centering,includefoot,foot=0.6in]{geometry}

\newtheorem{Lemma}{Lemma}[section]
\newtheorem{Corollary}[Lemma]{Corollary}

\newtheorem{Proposition}[Lemma]{Proposition}
\newtheorem{Remark}[Lemma]{Remark}
\newtheorem{Theorem}{Theorem}

\newenvironment{Proof}[1][.]%
 {\begin{trivlist}\item[]\textbf{Proof#1 }}%
 {\hspace*{\fill}$\rule{0.3\baselineskip}{0.35\baselineskip}$\end{trivlist}}

\newenvironment{Acknowledgment}%
 {\begin{trivlist}\item[]\textbf{Acknowledgments }}{\end{trivlist}}

\makeatletter\@addtoreset{equation}{section}\makeatother

\def\Re{\mathop\mathrm{Re}\nolimits}    
\def\Im{\mathop\mathrm{Im}\nolimits}    

\newcommand{\rmd}{\mathrm{d}}           
\newcommand{\rmi}{\mathrm{i}}           

\newcommand{\px}{\partial_{\xi}}

\newcommand{\si}{\sigma}

\newcommand{\be}{\begin{equation}}
\newcommand{\ee}{\end{equation}}
\newcommand{\bea}{\begin{eqnarray}}
\newcommand{\eea}{\end{eqnarray}}
\newcommand{\beas}{\begin{eqnarray*}}
\newcommand{\eeas}{\end{eqnarray*}}

\newcommand{\ba}{\begin{array}}
\newcommand{\ea}{\end{array}}

\newcommand{\bpm}{\begin{pmatrix}}
\newcommand{\epm}{\end{pmatrix}}

\begin{document}

\title{Rigorous justification of Taylor dispersion via center manifolds and hypocoercivity}

\author{Margaret Beck$^{\ast}$  \and Osman Chaudhary$^{\dagger}$  \and C. Eugene Wayne$^{\ddagger}$ }

\date{\today}


\maketitle


\begin{abstract}
Taylor diffusion (or dispersion) refers to a phenomenon discovered experimentally by Taylor in the 1950s where a solute dropped into a pipe with a background shear flow experiences diffusion at a rate proportional to $1/\nu$, which is much faster than what would be produced by the static fluid if its viscosity is $0 < \nu \ll 1$. This phenomenon is analyzed rigorously using the linear PDE governing the evolution of the solute. It is shown that the solution can be split into two pieces, an approximate solution and a remainder term. The approximate solution is governed by an infinite-dimensional system of ODEs that possesses a finite-dimensional center manifold, on which the dynamics correspond to diffusion at a rate proportional to $1/\nu$. The remainder term is shown to decay at a rate that is much faster than the leading order behavior of the approximate solution. This is proven using a spectral decomposition in Fourier space and a hypocoercive estimate to control the intermediate Fourier modes.
\end{abstract}


\section{Introduction}

Taylor dispersion is a phenomenon in fluid dynamics that was discovered in the 1950's by Geoffrey Taylor \cite{Taylor:1953,Taylor:1954}. The setting is a three dimensional pipe in which there is a background shear flow advecting the fluid down the length of the pipe, but where the rate of advection can vary as a function of the cross-sectional variables. It was observed by Taylor that, if a localized drop of dye was put into the pipe, then as expected it would be carried down the pipe by the shear flow and also diffuse due to the non-zero fluid viscosity. However, what was not expected was that the rate of diffusion experienced by the dye was not that of the fluid, say $\nu$, but instead a rate proportional to $1/\nu$, which is much larger if $0 < \nu \ll 1$. This phenomenon has been subsequently analyzed by many people, for example \cite{Aris:1956,Chatwin:1985,Mercer:1990}, but most of the work has been formal, based on asymptotic calculations. Our goal in this work is to rigorously analyze Taylor dispersion and provide a mathematical mechanism for its occurrence using center manifolds and Villani's theory of hypocoercivity \cite{Villani09}. We note there is another rigorous analysis of Taylor dispersion, \cite{BedrossianCoti-Zelati17}, that also uses hypoceorcivity in the proof. We will comment on the relationship between that and the present work at the end of this section. 

The PDE model of fluid flow in a pipe with a background shear flow is given by
\[
u_t = \nu \Delta u - V(y,z) u_x, \qquad x \in \mathbb{R}, \qquad (y,z) \in \Omega \subset \mathbb{R}^2. 
\]
The function $u: \mathbb{R} \times \Omega \times \mathbb{R}^+ \to \mathbb{R}$ represents the concentration of the solute, or dye, and the function $V: \Omega \to \mathbb{R}$ is a smooth background shear flow, which depends only on the cross-sectional variables $(y,z) \in \Omega$, where $\Omega$ is compact with smooth boundary. We assume Neumann boundary conditions, 
\[
\frac{\partial u}{\partial n}|_{\partial \Omega} = 0.
\]
For simplicity we assume the viscosity is a small, positive constant, $0 < \nu \ll 1$. To remove any effects of constant background advection caused by $V$, we define $\chi$ via 
\[
V(y,z) = A(1 + \chi(y,z)), \qquad A =\frac{1}{\mathrm{vol}(\Omega)} \int_{\Omega} V(y,z) \rmd y \rmd z, 
\]
and require that $\chi \in H^2(\Omega)$. Thus, $A$ is the average rate of advection in a cross section, and $\chi$ therefore has zero average advection in a cross section. We can then change variables using $x \to x + At$ to obtain 
\begin{equation}\label{E:orig-model}
u_t = \nu \Delta u - A\chi (y,z) u_x.
\end{equation}
It will be convenient to separate the effects of the cross-stream and longitudinal pipe variables. To that end, we will expand both $u$ and $\chi$ in terms of the eigenfunctions of the Laplacian $ \partial_y^2 + \partial_z^2$ acting on the compact domain $\Omega$. These eigenfunctions, which we denote by $\{\psi_n\}_{n=0}^\infty$, form an orthonormal basis for $L^2(\Omega)$ with $\psi_0 \equiv 1$, and we denote their corresponding eigenvalues by $\{-\mu_n\}_{n=0}^\infty$, which satisfy $0 = \mu_0 < \mu_1 \leq \mu_2 \leq \dots$ \cite[\S11.3]{Strauss08}. It will also be helpful to scale the longitudinal space variable $x$ and the time variable $t$ by $\nu$ via 
\begin{eqnarray} \label{E:nudep}
X = \nu x, \qquad T = \nu t.
\end{eqnarray}
This transforms (\ref{E:orig-model}) into 
\bea \label{E:scaled-model}
u_T = \nu^2 u_{XX} +  \Delta_{y,z} u - A\chi (y,z) u_X.
\eea
The main advantage of this change is that it helps us determine the dependence of the solutions on the viscosity parameter $\nu \ll 1$. This advantage will be made clear in Remarks \ref{rem:cmConst} and \ref{rmk:spec-nu}.
Inserting the expansions
\begin{equation}\label{E:initial-expansion}
u(X, y, z, T) = \sum_{n=0}^\infty u_n(X,T) \psi_n(y,z), \qquad \chi(y,z) = \sum_{n=0}^\infty \chi_n \psi_n(y,z),
\end{equation}
where
\[
u_n(X,T) = \int_\Omega u(X,y,z,T) \psi_n(y,z) \rmd y \rmd z, \qquad \chi_n = \int_\Omega \chi(y,z) \psi_n(y,z) \rmd y \rmd z,
\]
into equation \eqref{E:scaled-model} and noting that $\chi_0 = 0$ since it has zero average in $\Omega$, we obtain
\begin{eqnarray}
\partial_T u_0 &=& \nu^2 \partial_X^2 u_0  - A \sum_{m=1}^\infty \chi_m \partial_X u_m \label{E:u0} \\
\partial_T u_n &=& \nu^2\partial_X^2 u_n - \mu_n u_n - A \chi_n\partial_X u_0  - A \sum_{m = 1}^\infty \chi_{n, m} \partial_X u_m, \qquad n = 1, 2, \dots, \label{E:un}
\end{eqnarray}
where 
\[
\chi_{n, m} = \langle \psi_n, \chi \psi_m \rangle_{L^2(\Omega)}.
\]

In order to use invariant manifolds to study Taylor dispersion, we must deal with the fact that the Laplacian, $\partial_X^2$, on $\mathbb{R}$ has continuous spectrum consisting of $(-\infty, 0]$; in other words, there is no spectral gap. One way to overcome this is to use similarity variables, 
\[
\xi = \frac{X}{\sqrt{T+1}}, \qquad \tau = \log(T + 1),
\]
which exploit the space/time scaling inherent to the operator \cite{Wayne:1997}. (The use of $T+1$, rather than $T$, in the above definition is just for convenience, so that the change of variables is well-defined at $T = 0$.) We therefore further define new dependent variables $\{w_n\}_{n=0}^\infty$ via
\begin{eqnarray} \label{E:covScale}
u_0(X,T) && = \frac{1}{\sqrt{T+1}} w_0 \left( \frac{X}{\sqrt{T+1}},\log(T + 1)\right) \\
\qquad u_n(X,T) && = \frac{1}{T+1} w_n \left( \frac{X}{\sqrt{T+1}},\log(T + 1)\right), \qquad n = 1, 2, \dots.
\end{eqnarray}
Plugging this definition into \eqref{E:u0}-\eqref{E:un}, we obtain
\begin{eqnarray}
\partial_\tau w_0 &=& \mathcal{L} w_0  - A \sum_{m=1}^\infty \chi_m \partial_\xi w_m\label{E:w0}  \\ 
\partial_\tau w_n &=& \left(\mathcal{L} + \frac{1}{2}\right) w_n - e^{\tau/2} A \sum_{m = 1}^\infty \chi_{n, m} \partial_\xi w_m - e^\tau( \mu_n w_n + A \chi_n \partial_\xi w_0), \label{E:wn}
\end{eqnarray}
where
\begin{equation}\label{E:L}
\mathcal{L} = \nu^2 \partial_\xi^2 + \frac{1}{2} \partial_\xi(\xi \cdot) = \nu^2 \partial_\xi^2 + \frac{1}{2} \xi \partial_\xi + \frac{1}{2}
\end{equation}
is the Laplacian $\nu^2 \partial_X^2$ written in terms of the similarity variables. Note that the reason for the different powers of $(T+1)$ in front of $w_0$ and $w_n$ for $n \geq 1$ in (\ref{E:covScale}) is that equation (\ref{E:w0}) above becomes $\tau-$ independent. Continuing, we remark that the operator $\mathcal{L}$ was analyzed in detail in \cite{Gallay:2002}. Its properties are given in \S\ref{S:setupresults} below, but for the moment we just note that, on the space
\begin{equation}\label{E:L2m}
L^2(m) = \left\{ w \in L^2(\mathbb{R}): \int_\mathbb{R} (1 + \xi^2)^m |w(\xi)|^2 \rmd \xi < \infty \right\},
\end{equation}
the spectrum of $\mathcal{L}$ is composed of essential and discrete spectrum: 
\[
\sigma(\mathcal{L}) =  \{ \lambda \in \mathbb{C}: \mathrm{Re}(\lambda) \leq -(2m-1)/4\} \cup \{ \lambda = -k/2: k = 0, 1, 2, \dots\}.
\]
Thus, as the algebraic weight $m$ in the definition of the function space $L^2(m)$ increases, the essential spectrum is pushed further into the left half-plane, revealing more and more isolated eigenvalues at negative multiples of $1/2$. This suggests that we can construct a center-stable manifold (which we often refer to as a center manifold, for short) corresponding to those isolated eigenvalues, where the dimension of this manifold can be large if $m$ is sufficiently large. 

The utility of such a center manifold can be seen by considering the term $-e^\tau( \mu_n w_n + A \chi_n \partial_\xi w_0)$ in \eqref{E:wn}. As $\tau$ increases this term becomes large, which suggests that $w_n$ should evolve so that ultimately $\mu_n w_n + A \chi_n \partial_\xi w_0 = 0$. Hence, we expect that, for large times, 
\[
w_n \approx - \frac{A\chi_n}{\mu_n} \partial_\xi w_0 \qquad \Rightarrow \qquad \partial_\tau w_0 \approx \mathcal{L}_{td} w_0,
\]
where 
\begin{equation}\label{E:LT}
\mathcal{L}_{td}:= \left(\nu^2 + A^2 \|\chi\|_{\mu}^2\right) \partial_\xi^2 + \frac{1}{2} \partial_\xi(\xi \cdot), \qquad \|\chi\|_{\mu}^2 = \sum_{m} \frac{1}{\mu_m}\chi_m^2
\end{equation}
is again the Laplacian in similarity variables but now with Taylor diffusion coefficient
\begin{equation}\label{E:nuT}
\nu_{td} := \left(\nu^2 + A^2 \|\chi\|_{\mu}^2\right).
\end{equation}
Note that the spectrum of the operator does not depend on the viscosity, so $\sigma(\mathcal{L}) = \sigma(\mathcal{L}_{td})$. Thus, we expect that $\{w_n\}_{n=1}^\infty$ will rapidly converge to a manifold defined by $w_n = -(A\chi_n \partial_\xi w_0)/(\mu_n)$, and then for large times the dynamics of $w_0$ can be described by a center-stable manifold corresponding to the isolated eigenvalues of the operator $\mathcal{L}_{td}$. In terms of the original variables, this suggests that $\{u_n\}_{n=1}^\infty$ should become ``slaved" to the low mode $u_0$ exponentially fast, while the low mode $u_0$ should decay diffusively, but as if its diffusion coefficient is $\nu_{td} = \mathcal{O}(1)$ (instead of $\nu^2$), which, if we change back to the original $(x,t)$ variables, matches the experimental observations of Taylor and the formal calculations in \cite{Chatwin:1985}.  

There are several technical difficulties that must be overcome in order to make the above argument rigorous. First, in analyzing the dynamics of system \eqref{E:w0}-\eqref{E:wn} using the spectral structure of $\mathcal{L}_{td}$, it would be natural to expand each $w_n$, $n = 0, 1, \dots$, in terms of the eigenfunctions $\{\varphi_j^{td}(\xi)\}_{j=0}^N$ of $\mathcal{L}_{td}$, where $N = N(m)$ corresponds to the number of isolated eigenvalues, and hence the dimension of the center-stable manifold. In other words, we could write
\[
P_N w_n(\xi, \tau) = \sum_{j=0}^N \alpha_{j, n}(\tau) \varphi_j^{td}(\xi), \qquad w_n^\mathrm{s} = (1-P_N)w_n
\] 
for each $n$, where $w_n^\mathrm{s}$ is the component of the solution in the strong stable manifold, which we expect to decay rapidly. 
Although this is essentially what we will do, it turns out that it will be more convenient to prove the rapid decay of $w_n^\mathrm{s}$ in terms of the $(X, T)$ variables, by using the Fourier transform. 

The reason for this is that our center manifold argument will only show that the enhanced diffusion affects the first $N+1$ terms in the eigenfunction expansion $\{w_n\}$. This is sufficient for the physical realization of the phenomenon because the higher order terms, corresponding to $w_n^{s}$, will be shown to decay like $T^{-\mathcal{N}(N)}$ where $\mathcal{N}$ can be made large by choosing $N$, and hence also $m$, to be large, which is faster than the enhanced algebraic diffusive decay resulting from Taylor diffusion. 

To understand what $P_N w_n$ corresponds to in the physical $(x,t)$ variables, consider the following calculation. The eigenfunctions of $\mathcal{L}_{td}$ are given by
\begin{equation}\label{E:defvarphi}
\varphi_j^{td}(\xi) = \partial_\xi^j \varphi_0^{td}(\xi), \qquad \varphi_0^{td}(\xi) = \frac{1}{\sqrt{4\pi\nu_{td}}}e^{-\frac{\xi^2}{4\nu_{td}}}.
\end{equation}
If we assume that
\[
u(X,T) = \frac{1}{(1+T)^\gamma} w\left( \frac{X}{\sqrt{T+1}}, \log(T+1)\right), \qquad w(\xi, \tau) = \sum_{j=0}^N \alpha_j(\tau)\varphi_j^{td}(\xi), 
\]
which can represent either $w_0$ or $w_n$, $n \geq 1$, depending on the choice of $\gamma$, then 
\begin{eqnarray*}
\hat u(\kappa,T) &=& \int e^{\rmi \kappa X} u(X,T) \rmd X \\
&=& \sum_{j=0}^N \frac{1}{\sqrt{4\pi\nu_{td}}}(1+T)^{j/2-\gamma} (-\rmi \kappa)^j \alpha_j(\log(T+1))\int e^{\rmi \kappa X} e^{-\frac{X^2}{4\nu_{td}(T+1)}} \rmd X \\
&=& \sum_{j=0}^N (1+T)^{j/2+1/2-\gamma} (-\rmi \kappa)^j \alpha_j(\log(T+1)) e^{-\kappa^2\nu_{td}(T+1)}.
\end{eqnarray*}
This implies that
\[
\hat u(0, T) = (1+T)^{1/2-\gamma}\alpha_0(\log(T+1)), \qquad \hat \partial_{\kappa} u(0, T) = (-\rmi) (1+T)^{1-\gamma}\alpha_1(\log(T+1)), \qquad \dots 
\]
which, combined with the Taylor expansion
\[
\hat u(\kappa, T) = \hat u(0, T) + \partial_\kappa \hat u(0, T) \kappa + \frac{1}{2} \partial_\kappa^2 \hat u(0, T) \kappa^2 + \dots 
\]
means that the behavior of $P_N w$ tells us about the behavior of $\hat u(\kappa, T)$ for $\kappa$ near zero. In other words, $P_N w$ represents both the behavior of the ``low modes" of $w(\xi, \tau)$, where ``low modes" refers to the leading eigenfunctions of $\mathcal{L}_{td}$, and the behavior of the ``low modes" of $\hat u(\kappa,T)$, where now ``low modes" refers to values of the Fourier variable $\kappa$ near zero. This relationship between Taylor dispersion, the behavior of the Fourier transform of the solution at small wave numbers, and the center-manifold theorem was also discussed by Mercer and Roberts in \cite{Mercer:1990}.

We will refer to $w^\mathrm{s} = (1-P_N)w$ as the remainder, or error, term. In terms of $\mathcal{L}_{td}$ it corresponds to the behavior due to the essential spectrum and the discrete spectrum that is sufficiently far from the imaginary axis. To prove that the remainder term decays rapidly, it will be convenient to work in terms of the Fourier variables associated with physical $(X,T)$ space, rather than system \eqref{E:w0}-\eqref{E:wn}. This will lead to a linear, nonautonomous equation governing the behavior of the remainder term of the form $\hat U_T = \mathcal{B}(\kappa) \hat U + \hat F(\kappa, T)$. We can then consider three regimes: a small wavenumber regime defined by $|\kappa| \leq \kappa_0$, an intermediate one defined by $\kappa_0 \leq |\kappa| \leq \kappa_1/\nu$, and a large one defined by $|\kappa| \geq \kappa_1/\nu$. In the large regime, the solution decays exponentially due to the usual (non-Taylor) diffusive estimate $e^{\nu^2 \partial_X^2 T} \sim e^{-\nu^2 \kappa^2 T} \leq e^{- \kappa_1^2 T}$. In the intermediate regime this naive estimate is not quite strong enough, because it only implies $ e^{\nu^2 \partial_X^2 T} \sim e^{-\nu^2 \kappa^2 T} \leq \sim e^{-\nu^2 \kappa_0^2 T}$, which is quite weak for $0 < \nu \ll 1$. To improve it, we will apply a hypocoercivity argument \cite{Villani09} to show that in this region we also have decay like $e^{-M T}$ for some $M > 0$. For the low wavenumbers, we will decompose the remainder term into a piece corresponding to the leading eigenvalue $\lambda_0(\kappa)$ of $\mathcal{B}(\kappa)$, which is parabolic with $\lambda_0(0) = 0$, and a piece corresponding to the rest of the spectrum of $\mathcal{B}(\kappa)$. The latter will decay exponentially fast because $\mathcal{B}(k)$ has a spectral gap for each fixed $k$. The former will be shown to decay algebraically with the rate $T^{-\mathcal{N}(N)}$, because we have already removed the leading order behavior via the term $P_N w_n$. 

Our analysis will be divided into the following steps. In \S\ref{S:setupresults} we will more precisely set-up our problem and carefully state the main results. In \S\ref{S:lowmodes} we will use the similarity variables and a center-stable manifold to prove that the low modes, corresponding to $P_Nw_n$, experience enhanced Taylor diffusion. Finally, in \S\ref{S:remainder} we will use a spectral decomposition and hypocoercivity to show that the remainder term decays rapidly, thus allowing for the Taylor diffusion to be physically observable. 

Before carrying this out, we comment on other related rigorous work on Taylor diffusion. In \cite{BeckChaudharyWayne15} we analyzed a model of system \eqref{E:u0}-\eqref{E:un} consisting of only two equations, one corresponding to $u_0$ and one modeling all of the $u_n$ for $n \geq 1$, and carried out a similar analysis there. This allowed us to focus on the main ideas of the argument: that the Taylor diffusion is really only affecting the low modes, with the remainder term decaying rapidly. However, in that work, because of the simple form of the system, one could see directly that the remainder term decayed rapidly and the hypocoercivity argument we use here in \S\ref{S:remainder} was not necessary. Moreover, the center manifold argument, which was used to justify the enhanced diffusion, was constructed for a finite-dimensional ODE. Here, the center manifold argument in \S\ref{S:lowmodes} will need to be carried out for an infinite-dimensional ODE.

Also, in \cite{BedrossianCoti-Zelati17} an equation very similar to \eqref{E:orig-model} was analyzed, also using hypocoercivity. However, there Villani's framework was applied directly to the PDE \eqref{E:orig-model}, whereas our hypocoercivitiy argument is applied in Fourier space. This allows us to avoid any assumptions on the critical points of the shear flow $\chi$, which play an important role in the argument in \cite{BedrossianCoti-Zelati17}. Moreover, since $X \in \mathbb{R}$, we need to work in Fourier space with all $|\kappa| \geq 0$. The setting in \cite{BedrossianCoti-Zelati17} is for a bounded $X$ domain, which effectively means $|\kappa| \geq 1$. This changes the nature of the resulting decay and the regions in which the enhanced diffusion is obtained.


\section{Set-up and statement of main results} \label{S:setupresults}

The main result that we will prove is the following. Theorem \ref{thm:main1}(i) will be proven in \S\ref{S:lowmodes} and Theorem \ref{thm:main1}(ii) will be proven in \S\ref{S:remainder}. In the statement of the Theorem we use the following notation for the space in which the initial data must lie:
\[
L^2(N+1)\times L^2(\Omega) = \left\{ u \in L^2(\mathbb{R} \times \Omega): \int_\mathbb{R} \int_\Omega (1+X^2)^{N+1}|u(X, y, z)|^2 \rmd X \rmd y \rmd z =: \|u\|^2_{L^2(N+1)\times L^2(\Omega)} < \infty \right\}.
\]

\begin{Theorem} \label{thm:main1}  Given any $N > 0$, if $u(\cdot,0) \in L^2(N+1)\times L^2(\Omega)$, then there exist constants \newline $C_j= C_j(\| u( \cdot, 0) \|_{L^2(N+1)\times L^2(\Omega)})$, $j = 1,2$, that are independent of $\nu$ and a decomposition of the corresponding solution of \eqref{E:scaled-model} of the form
\[
u(X,y,z,T) = u_\mathrm{app}(X,y,z,T) + u_\mathrm{rem}(X,y,z,T),
\]
where $u_\mathrm{app}(X,y,z,T)$ and $u_\mathrm{rem}(X,y,z,T)$ are defined in \eqref{E:defuapp}-\eqref{E:defurem}, that satisfies the following. 
\begin{enumerate}
\item There exists an infinite-dimensional system of ordinary differential equations that govern the behavior of $u_\mathrm{app}$. Moreover, this system of ODEs possesses a finite dimensional center manifold that is globally attracting at a rate that is exponential in $T$, $e^{-\eta T}$ for some $\eta$ independent of $\nu$, and on which the dynamics correspond to enhanced diffusion with viscosity $\nu_{td}$, defined in \eqref{E:nuT}. In other words, 
\[
\left\|u_\mathrm{app}(X,y,z,T) - \frac{C_1}{\sqrt{4\pi\nu_{td}(T+1)}}e^{-\frac{X^2}{4\nu_{td}(T+1)}}\right\|_{L^2(\mathbb{R}\times \Omega)} \leq \frac{C_2}{ (1+T)^{3/4}},
\]
{\color{black} The constant $C_1$ is given explicitly by
\[
C_1 = \int_\mathbb{R}\int_\Omega u(X, y, z, 0)\rmd X \rmd y \rmd z.
\]
} 
\item The remainder term satisfies
\[
\|u_\mathrm{rem}(\cdot,T)\|_{L^2(\mathbb{R}\times \Omega)} \leq \frac{C_2}{(1+T)^{\frac{N}{6} + \frac{1}{12}}}.
\]
\end{enumerate}
\end{Theorem}

If we translate these results back to our original, unscaled time and space variables and choose
$N \ge 4$, so that $\frac{N}{6}+ \frac{1}{12} \ge \frac{3}{4}$,  we see that we obtain immediately:

{\color{black}
\begin{Corollary} Given any initial condtion $u(\cdot, 0) \in L^2(N+1)\times L^2(\Omega)$, there exist constants $\tilde C_j= \tilde C_j(\| u( \cdot, 0) \|_{L^2(N+1)\times L^2(\Omega)})$, $j = 1,2$, such that the solution of \eqref{E:orig-model} satisfies
$$
\left\| u(x,y,z,t) - \frac{\tilde C_1}{\sqrt{4\pi (\nu + A^2 \| \chi \|^2_{\mu}/\nu ) (t + 1/\nu)}} e^{- \frac{x^2}{4(\nu + A^2 \| \chi \|^2_{\mu}/\nu ) (t + 1/\nu)}}
\right \|_{L^2(\mathbb{R}\times \Omega)}
\le \frac{\tilde C_2}{(1+ \nu t)^{3/4} }\ .
$$
The constant $\tilde C_1$ is given explicitly by
\[
\tilde C_1 = \int_\mathbb{R}\int_\Omega u(x, y, z, 0)\rmd x \rmd y \rmd z.
\]
\end{Corollary}
}

\begin{Remark}
Note that the leading order term in the asymptotics identified by this Corollary corresponds to a solution
of the diffusion equation with diffusion coefficient $(\nu + A^2 \| \chi \|_{\mu}^2/\nu )$ which is precisely
the asymptotic behavior derived non-rigorously in \cite{Smith:1987}. (In particular, see (2.17) for the calculation of the shear diffusion coefficient.) {\color{black} We note that the constant $\tilde C_2$ appearing in the Corollary can be related to the constant $C_2$ appearing in the Theorem by undoing the change of variables $X = \nu x$.}
\end{Remark}

\begin{Remark}  
As we discuss later in Section 3, we actually derive not just the leading order term in the asymptotics but higher terms as well - in principle, terms of arbitrary order, if  the initial condition $u_0$ decays
sufficiently rapidly as $|x| \to \infty$.  The higher order terms in the asymptotics are expressed in terms
of the eigenfunctions of the operator $\mathcal{L}_{td}$. {\color{black} See Remark \ref{rem:asym} for further details.}
\end{Remark}

To prove these results, we will use the following facts about the operator $\mathcal{L}_{td}$, which is just the Laplacian written in terms of similarity variables. Recall from \eqref{E:LT} that
\begin{equation*}
\mathcal{L}_{td}\varphi = \nu_{td} \partial_{\xi}^2 \varphi + \frac{1}{2} \partial_{\xi} (\xi \varphi).
\end{equation*}
We state the following results for viscosity $\nu_{td}$, but the results are true with $\nu_{td}$ replaced by any other positive number. This operator has been analyzed in \cite{Gallay:2002}, and in the weighted
Hilbert space $L^2(m)$ defined in \eqref{E:L2m} one finds
\begin{equation*}
\sigma(\mathcal{L}_{td}) = \left\{ \lambda \in \mathbb{C}: {\Re}(\lambda) \le - \frac{(2m-1)}{4}  \right\} \cup \left\{ -\frac{k}{2} ~|~ 
k \in {\mathbb{N}}  \right\} \ .
\end{equation*}
Furthermore, the eigenfunctions corresponding to the isolated eigenvalues $\lambda_k = -k/2$ are given by
the Hermite functions
\begin{equation*}
\varphi_0^{td}(\xi) = \frac{1}{\sqrt{4 \pi \nu_{td} } } e^{-\frac{\xi^2}{4 \nu_{td}}}, \qquad   \varphi_k^{td}(\xi) =  \partial_{\xi}^k \varphi_0^{td}(\xi). 
\end{equation*}
The corresponding adjoint eigenfunctions are given by the Hermite polynomials
\begin{equation}\label{E:hkt}
H_k^{td}(\xi) = \frac{ 2^k \nu_{td}^k}{k ! } e^{\frac{\xi^2}{4 \nu_{td}}}\partial_{\xi}^k e^{-\frac{\xi^2}{4 \nu_{td}}}\ .
\end{equation}
Note that we have the orthogonality relationship
\[
\langle H_k^{td}, \varphi_j^{td} \rangle_{L^2(\mathbb{R})} = \delta_{jk} = \begin{cases} 1 & \mbox{ if } j = k \\ 0 & \mbox{ if } j \neq k,\end{cases}
\]
which can be used to define spectral projections.

\begin{Remark}  The expressions in \cite{Gallay:2002} for $\varphi_k^{td}$ and $H_k^{td}$ are derived in the case when $\nu_{td} = 1$.  The expressions given here follow easily by the change of variables 
$\xi \to \xi/\sqrt{\nu_{td}}$.  
\end{Remark}


\subsection{Preparation of the equations}\label{S:eqnprep}

To emphasize the expected role of the enhanced diffusion, we rewrite \eqref{E:w0}-\eqref{E:wn} as
\begin{eqnarray}
\partial_\tau w_0 &=& \mathcal{L}_{td} w_0  - D_{td}\partial_\xi^2 w_0 - A \sum_{m=1}^\infty \chi_m \partial_\xi w_m \label{E:w0T2} \\
\partial_\tau w_n &=& \left(\mathcal{L}_{td} + \frac{1}{2}\right) w_n - D_{td}\partial_\xi^2 w_n - e^{\tau/2} A \sum_{m = 1}^\infty \chi_{n, m} \partial_\xi w_m - e^\tau( \mu_n w_n + A \chi_n \partial_\xi w_0), \label{E:wnT}
\end{eqnarray}
where 
\[
D_{td} := A^2 \|\chi\|_\mu^2
\]
and $\mathcal{L}_{td}$ is defined in equation \eqref{E:LT}. As described above, asymptotically we expect $w_n = -(A\chi_n \partial_\xi w_0)/(\mu_n)$, which is a perfect derivative. To exploit this, we wish to effectively integrate the $w_n$ equation. Naively, this could be done by defining $\{V_n\}_{n=1}^\infty$ via $\partial_\xi V_n = w_n$. In order to obtain decay of $V_n$ as $|\xi| \to \infty$, we would then need to assume that $\int w_n = 0$. To avoid this additional assumption, we instead define $\{V_n\}_{n=1}^\infty$ via
\begin{equation}\label{E:defvnT}
w_n(\xi, \tau) = \gamma_n(\tau) \varphi_0^{td}(\xi) + V_n(\xi, \tau), \qquad \gamma_n(\tau) = \int_\mathbb{R} w_n(\xi, \tau)\rmd \xi,
\end{equation}
where $\varphi_0^{td}$ is the eigenfunction of $\mathcal{L}_{td}$ defined in \eqref{E:defvarphi} associated with the zero eigenvalue. 
Note that this implies
\[
\gamma_n(\tau) = \langle w_n(\tau), H_0^{td} \rangle_{L^2} = \int_\mathbb{R} w_n(\xi, \tau) \rmd \xi
\]
and that $\gamma_n(\tau)$ is bounded for each $\tau$ such that $w_n(\tau) \in L^2(m)$, with $m > 1/2$, because
\[
|\gamma_n(\tau)| \leq \int_\mathbb{R} \frac{1}{(1+\xi^2)^\frac{m}{2}}(1+\xi^2)^\frac{m}{2} |w(\xi, \tau)| \rmd \xi \leq \left( \int_\mathbb{R} \frac{1}{(1+\xi^2)^m} \rmd \xi \right)^{1/2} \|w_n(\tau)\|_{L^2(m)} \leq C(m) \|w_n(\tau)\|_{L^2(m)}.
\]
Since $\int \varphi_0^{td} = 1$, we see that $\int V_n = 0$. Inserting \eqref{E:defvnT} into \eqref{E:wnT}, we find
\begin{eqnarray}
\dot \gamma_n \varphi_0^{td} + \partial_\tau V_n &=& \frac{1}{2} \gamma_n \varphi_0^{td} + \left(\mathcal{L}_{td} + \frac{1}{2}\right) V_n 
- D_{td}\left( \gamma_n \varphi_2^{td} + \partial_\xi^2 V_n\right) \label{E:vnTtemp-2} \\
&& \qquad \qquad - e^{\tau/2} A \sum_{m = 1}^\infty \chi_{n, m} \partial_\xi (\gamma_m \varphi_0^{td} + V_m) 
- e^\tau(  \mu_n V_n + A \chi_n \partial_\xi w_0) - e^\tau \mu_n \gamma_n \varphi_0^{td}. \nonumber
\end{eqnarray}
Integrating over $\mathbb{R}$ and using the fact that $\varphi_k^{td}, w_0 \to 0$ as $|\xi| \to \infty$, we find
\[
\dot \gamma_n = \left( \frac{1}{2} - e^\tau \mu_n\right) \gamma_n, 
\]
which implies that 
\begin{equation}\label{E:gamman}
\gamma_n(\tau) = \gamma_n(0) e^{\frac{\tau}{2} - \mu_n(e^\tau - 1)}.
\end{equation}
With this information, in  \eqref{E:vnTtemp-2} we can cancel all the terms involving $\gamma_n$ alone, use the fact that $\int V_n = 0$ to define $v_n$ via $\partial_\xi v_n = V_n$, and obtain from \eqref{E:w0T2}-\eqref{E:wnT}
\begin{eqnarray}
\partial_\tau w_0 &=& \mathcal{L}_{td} w_0  - D_{td}\partial_\xi^2 w_0  - A  \sum_{m=1}^\infty \chi_m \partial_\xi^2 v_m - A e^\frac{\tau}{2} \varphi_1^{td}   \sum_{m=1}^\infty \chi_m \gamma_m(0) e^{- \mu_m(e^\tau - 1)}  \label{E:w0T3} \\
\partial_\tau v_n &=& \mathcal{L}_{td} v_n - D_{td}\partial_\xi^2 v_n - e^{\tau/2} A \sum_{m = 1}^\infty \chi_{n, m} \partial_\xi v_m - e^\tau( \mu_n v_n + A \chi_n w_0) \nonumber \\
&&\qquad \qquad - e^{\tau} A \varphi_0^{td}  \sum_{m = 1}^\infty \chi_{n, m} \gamma_m(0) e^{ - \mu_m(e^\tau - 1)} -  D_{td} \gamma_n(0) e^{\frac{\tau}{2}} e^{-\mu_n(e^\tau-1)} \varphi_1^{td}. \label{E:vnT2}
\end{eqnarray}


\subsection{Separation into low modes and the remainder term}

In order to analyze the behavior of solutions to system \eqref{E:w0T3}-\eqref{E:vnT2}, we define
\begin{eqnarray}
w_0(\xi, \tau) = \sum_{k=0}^N \alpha_k(\tau) \varphi_k^{td}(\xi) + w_0^s(\xi, \tau) \nonumber \\
v_n(\xi, \tau) = \sum_{k=0}^N \beta_k^n(\tau) \varphi_k^{td}(\xi) + v_n^s(\xi, \tau), \label{E:deflowhigh}
\end{eqnarray}
where $\{\varphi_k^{td}\}_{k=0}^N$ are the first $N+1$ eigenfunctions associated with $\mathcal{L}_{td}$ and 
\[
\alpha_k(\tau) = \langle w_0(\xi, \tau), H_k^{td}(\xi) \rangle_{L^2(\mathbb{R})}, \qquad \beta_k^n(\tau) = \langle v_n(\xi, \tau), H_k^{td}(\xi) \rangle_{L^2(\mathbb{R})},
\]
are the spectral projections onto those eigenmodes defined via the corresponding adjoint eigenfunctions $H_k^{td}$. See \eqref{E:hkt}.      Recalling that $\partial_\xi \varphi_k^{td} = \varphi_{k+1}^{td}$ and $\mathcal{L}_{td} \varphi_k^{td} = - (k/2) \varphi_k^{td}$, inserting the above expressions into 
\eqref{E:w0T3}-\eqref{E:vnT2} and taking the inner product of the result with $H_k^{td}$ gives the following infinite-dimensional system of ODEs for the evolution of $\{\alpha_k\}_{k=0}^N$ and $\{\beta_k^n\}_{k=0}^N$, $n \geq 1$:
\begin{eqnarray}
\dot \alpha_0 &=& 0 \nonumber \\
\dot \alpha_1 &=& -\frac{1}{2} \alpha_1 - A e^\frac{\tau}{2} \sum_{m=1}^\infty \chi_m \gamma_m(0) e^{-\mu_m(e^\tau - 1)} \nonumber \\
\dot \alpha_k &=& - \frac{k}{2} \alpha_k - D_{td}\alpha_{k-2} - A \sum_{m=1}^\infty \chi_m \beta_{k-2}^m \qquad \qquad  2 \leq k \leq N \nonumber \\
\dot \beta_0^n &=& -e^\tau(\mu_n \beta_0^n + A\chi_n \alpha_0) -e^\tau A \sum_{m=1}^\infty \chi_{n,m}\gamma_m(0) e^{-\mu_m(e^\tau-1)} \label{E:alphabeta}\\
\dot \beta_1^n &=& - \frac{1}{2}\beta_1^n - e^\tau( \mu_n \beta_1^n + A\chi_n \alpha_1) - e^\frac{\tau}{2} A \sum_{m=1}^\infty \chi_{n,m} \beta_0^m  - D_{td} \gamma_n(0)e^{\frac{\tau}{2}}e^{-\mu_n(e^\tau-1)}  \nonumber \\
\dot \beta_k^n &=& - \frac{k}{2} \beta_k^n - e^\tau( \mu_n \beta_k^n + A\chi_n \alpha_k) - D_{td} \beta_{k-2}^n  - e^\frac{\tau}{2} A \sum_{m=1}^\infty \chi_{n,m} \beta_{k-1}^m \qquad \qquad  2 \leq k \leq N. \nonumber 
\end{eqnarray}
Note that we have used the following facts. First, $\langle H_k^{td}, w_0^s \rangle_{L^2} = 0$, which follows by construction. This implies that $\langle H_k^{td}, \mathcal{L}_{td} w_0^s \rangle_{L^2} = \langle -(k/2) H_k^{td}, w_0^s \rangle_{L^2} = 0$. One can also 
check that
\[
\mathcal{L}_{td}^*(\partial_\xi H_k^{td}) = -\frac{(k-1)}{2} \partial_\xi H_k^{td} \qquad \Rightarrow \qquad H_{k-1}^{td} = - \partial_\xi H_k^{td}\,
\]
which implies that $\langle H_k^{td}, \partial_\xi^2 w_0^s \rangle = 0$. Similar results hold for $v_n^s$. 

The key aspect of \eqref{E:alphabeta} is that, because of the structure of \eqref{E:w0T3}-\eqref{E:vnT2}, the dynamics of $\{\alpha_k\}_{k=0}^N$ and $\{\beta_k^n\}_{k=0}^N$ do not depend on the remainder terms $w_0^\mathrm{s}$ or $v_n^\mathrm{s}$. Therefore, the behavior of these low modes can be analyzed without any a priori knowledge of the remainder terms. The structure of the above system suggests that, with the exception of $\alpha_0$, everything should decay exponentially fast in $\tau$, which corresponds to algebraic decay in $t$. Moreover, the leading order behavior will be governed by $\alpha_0$.


\subsection{Definition of $u_\mathrm{app}$ and $u_\mathrm{rem}$}

We now relate the decomposition in \eqref{E:deflowhigh} back to the solution $u(X,y,z,T)$ of the original equation \eqref{E:scaled-model}. We define $u_\mathrm{app}$ in terms of the low modes and $u_\mathrm{rem}$ in terms of the functions $w_0^s$ and $v_n^s$. To do so we need to convert back to the $(X,T)$ variables and take into account the decomposition in \eqref{E:initial-expansion} and the change of variables in \S\ref{S:eqnprep}. In particular, we have
\begin{eqnarray*}
u(X,y,z,T) &=& \sum_{n=0}^\infty u_n(X,T) \psi_n(y,z) \nonumber \\
u_0(X,T) &=& \frac{1}{\sqrt{1+T}}w_0(\xi, \tau)  \\
u_n(X,T) &=& \frac{1}{(T+1)} \left[ \gamma_n(\tau) \varphi_0^{td}(\xi) + \partial_\xi v_n(\xi, \tau)\right], \qquad n \geq 1. \nonumber
\end{eqnarray*}
Using \eqref{E:deflowhigh}, we find
\begin{eqnarray*}
u_0(X,T) &=& \frac{1}{\sqrt{1+T}} \sum_{k=0}^N \alpha_k(\tau) \varphi_k^{td}(\xi) + \frac{1}{\sqrt{T+1}} w_0^s(\xi, \tau) \\
u_n(X,T) &=& \frac{1}{(T+1)}\left[ \gamma_n(\tau) \varphi_0^{td}(\xi) + \sum_{k=0}^N \beta_k^n(\tau) \varphi_{k+1}^{td}(\xi) \right] + \frac{1}{(T+1)} \partial_\xi v_n^s(\xi, \tau).
\end{eqnarray*}
We now define
\begin{eqnarray}
u_\mathrm{app}(X,y,z,T) &=&  \frac{\psi_0(y,z)}{\sqrt{1+T}} \sum_{k=0}^N \alpha_k[\log(T+1)] \varphi_k^{td}\left(\frac{X}{\sqrt{T+1}}\right) \label{E:defuapp} \\
&& \qquad + \sum_{n=1}^\infty \frac{\psi_n(y,z)}{(T+1)} \left[ \gamma_n[\log(T+1)] \varphi_0^{td}\left(\frac{X}{\sqrt{T+1}}\right)  + \sum_{k=0}^N \beta_k^n[\log(T+1)] \varphi_{k+1}^{td}\left(\frac{X}{\sqrt{T+1}}\right)  \right], \nonumber
\end{eqnarray}
and
\begin{eqnarray}
u_\mathrm{rem}(X,y,z,T) &=&\frac{\psi_0(y,z)}{\sqrt{1+T}} w_0^s\left( \frac{X}{\sqrt{T+1}}, \log(T+1)\right) 
 + \sum_{n=1}^\infty \frac{\psi_n(y,z)}{(T+1)} \partial_\xi v_n^s \left( \frac{X}{\sqrt{T+1}}, \log(T+1)\right). \label{E:defurem} 
\end{eqnarray}
The behavior of $u_\mathrm{app}$, as stated in Theorem \ref{thm:main1}(i), will be determined in \S\ref{S:lowmodes}, and the behavior of $u_\mathrm{rem}$, as stated in Theorem \ref{thm:main1}(ii), will be determined in \S\ref{S:remainder}.

\begin{Remark}\label{rem:asym} {\color{black} Equation \eqref{E:defuapp} provides a way to compute higher order asymptotics of the solution. The leading order term, which appears in Theorem \ref{thm:main1}(i), corresponds only to the $\alpha_0$ term in \eqref{E:defuapp}. The functions $\alpha_j$, for $j = 1, \dots N$, as well as $\{\gamma_n\}$ and $\{\beta_k^n\}$, determine the higher order asymptotics. Indeed, one of the advantages of the center-manifold approach is that, in principle, we can compute the asymptotic behavior of the solution to any order.  More precisely, in \cite{Chaudhary17} it is proven that for any fixed inverse power of $t$, one can compute the behavior of the solution up to corrections of that order in $t$, solely in terms of the behavior of the solution restricted to the center manifold, which is given by the functions $\{\alpha_j\}$, $\{\gamma_n\}$, and $\{\beta_k^n\}$. Furthermore, because the formula for the center manifold, given in the proof of Proposition \ref{prop:cm}, is explicit, these functions could in principle also be computed explicitly.}
\end{Remark}


\section{Taylor dispersion for the approximate solution via a center manfold}\label{S:lowmodes}

The main goal of this section is to prove Theorem \ref{thm:main1}(i). This will essentially be done via Proposition \ref{prop:cm}, and it will be explained in \S\ref{S:thm1i} how its proof follows from that Proposition. 


\subsection{Asymptotic behavior of the low modes via a center-stable manifold}\label{S:low-cm}

Consider system \eqref{E:alphabeta}. To construct its center manifold, we start by performing some changes of variables. Recall from the formal analysis that, in long time limit, we expect $\mu_n w_n + A \chi_n \px w_0 = 0$. In system \eqref{E:alphabeta}, this results from the term $e^\tau(\mu_n \beta_k^n + A \chi_n \alpha_k)$. Therefore, we will diagonalize the system so that, in terms of new variables $(a_k, b_k^n)$, the set $\{ \mu_n \beta_k^n + A \chi_n \alpha_k\} = 0$ corresponds to the set $\{ b_k^n = 0\}$. We define 
\begin{equation} \label{eq:diagchange}
a_k  = \alpha_k, \qquad b_k^n  = \beta_k^n + \frac{A \chi_n}{\mu_n} \alpha_k
\end{equation}
and obtain
\begin{eqnarray*} 
\dot a_0 &=& 0 \\
\dot a_1 &=& -\frac{1}{2} a_1 - A e^\frac{\tau}{2} \sum_{m=1}^\infty \chi_m \gamma_m(0) e^{-\mu_m(e^\tau - 1)} \\
\dot{a}_k  &=& -\frac{k}{2} a_k - A \sum_{m=1}^{\infty} \chi_m b_{k-2}^m \qquad \qquad 2 \leq k \leq N \\
\dot{b}_0^n &=& -e^\tau \mu_n b_0^n - e^\tau A \sum_{m=1}^\infty \chi_{n,m} \gamma_m(0) e^{-\mu_m(e^\tau-1)} \\
\dot b_1^n &=& - \left(\frac{1}{2} + e^\tau \mu_n\right) b_1^n  - e^\frac{\tau}{2} A \sum_{m=1}^\infty \chi_{n,m} \left[ b_0^m - \frac{A \chi_m}{\mu_m} a_0\right]  - D_{td} \gamma_n(0)e^{\frac{\tau}{2}}e^{-\mu_n(e^\tau-1)} \\
&& \qquad \qquad  - \frac{A^2 \chi_n}{\mu_n} e^\frac{\tau}{2} \sum_{m=1}^\infty \chi_m \gamma_m(0) e^{-\mu_m(e^\tau-1)} \\
\dot{b}_k^n  &=& -\left( \frac{k}{2} + e^{\tau} \mu_n\right) b_k^n - D_{td} \left( b_{k-2}^n - \frac{A \chi_n}{\mu_n} a_{k-2} \right)  - \frac{A^2 \chi_n}{\mu_n} \sum_{m=1}^{\infty} \chi_m b_{k-2}^m  \\
&& \qquad \qquad - e^{\tau/2}A \sum_{m=1}^{\infty} \chi_{n,m} \left( b_{k-1}^m -\frac{A \chi_m}{\mu_m} a_{k-1} \right) \qquad \qquad 2 \leq k \leq N,
\end{eqnarray*}
where $n \geq 1$.
\begin{Remark} {\color{black} The equation for $\dot b_k^n$ follows from the fact that
\[
-D_{td} \left( \frac{A\chi_n}{\mu_n} \alpha_{k-2} + \beta_{k-2}^n \right) - \frac{A^2 \chi_n}{\mu_n} \sum_{m=1}^\infty \chi_m \beta_{k-2}^m =  -D_{td}\left( b_{k-2}^n - \frac{A \chi_n}{\mu_n} a_{k-2} \right)  - \frac{A^2 \chi_n}{\mu_n} \sum_{m=1}^{\infty} \chi_m b_{k-2}^m.
\]}
\end{Remark}
This system is non-autonomous, which makes it difficult to construct a center manifold. To overcome this, we first undo the change of variables in time using $\tau = \log(1+T)$ and define $\si = (1+T)^{-1/2}$. Denoting $d/dT = (\cdot)'$, we obtain
\begin{eqnarray} 
a_0' &=& 0 \nonumber \\
a_1' &=& -\frac{1}{2} \sigma^2 a_1 - A\sigma \sum_{m=1}^\infty \chi_m \gamma_m(0) e^{-\mu_mT} \nonumber \\
a_k'  &=& \sigma^2\left( -\frac{k}{2}a_k - A \sum_{m=1}^{\infty} \chi_m b_{k-2}^m \right) \qquad \qquad 2 \leq k \leq N \label{eq:diag}\\
{b_0^n}' &=& -\mu_n b_0^n - A \sum_{m=1}^\infty \chi_{n,m} \gamma_m(0) e^{-\mu_mT} \nonumber \\
{b_1^n}' &=& -\left(\frac{1}{2}\sigma^2 + \mu_n \right) b_1^n - A \sigma \sum_{m=1}^{\infty} \chi_{n,m} \left( b_{0}^m -\frac{A \chi_m}{\mu_m} a_{0}  \right) - D_{td} \sigma \gamma_n(0) e^{-\mu_n T} - \sigma \frac{A^2\chi_n}{\mu_n} \sum_{m=1}^\infty \chi_m \gamma_m(0) e^{-\mu_m T} \nonumber \\
{b_k^n}'  &=& -\left(\frac{k}{2}\sigma^2 +  \mu_n\right) b_k^n - D_{td} \sigma^2 \left( b_{k-2}^n - \frac{A \chi_n}{\mu_n} a_{k-2} \right)  - \frac{A^2 \chi_n}{\mu_n} \sigma^2 \sum_{m=1}^{\infty} \chi_m b_{k-2}^m \nonumber \\
&& \qquad \qquad - \sigma A \sum_{m=1}^{\infty} \chi_{n,m} \left( b_{k-1}^m -\frac{A \chi_m}{\mu_m} a_{k-1} \right) \qquad \qquad 2 \leq k \leq N \nonumber \\
\si'  &=& -\frac{1}{2} \si^3, \nonumber
\end{eqnarray}
where $n \geq 1$. Note that, except for the terms involving $\gamma_n(0)$, which are decaying exponentially fast in $T$, this system is autonomous (but nonlinear), due to our definition of $\sigma$.

It is now convenient to define more compact notation. To that end, we write
\[
b_k=(b_k^1, b_k^2, b_k^3, \ldots), \qquad \check{\chi}=(\chi_1, \chi_2, \chi_3 \ldots), \qquad \gamma = (\gamma_1, \gamma_2, \gamma_3, \dots), 
\]
where $\check\chi$ is a constant, $\gamma = \gamma(T)$ with $\gamma_n(T) = \gamma_n(0)e^{- \mu_n T}$, and $b_k = b_k(T)$, and $n \geq 1$. We also define operators on $\ell^2$ via
\[
(\tilde{\chi} * Y)_n = \sum_m \chi_{n,m} Y_m, \qquad (\Upsilon Y)_n = \mu_n Y_n.
\]
Throughout the following estimates we will use the following Lemma, which says that $\tilde{\chi}$ and $\Upsilon^{-1}$ are bounded operators.

\begin{Lemma}  The operators $\tilde{\chi}$ and $\Upsilon^{-1}$ are bounded operators on $\ell^2$.
\end{Lemma}

\begin{Proof}
The bound on $\Upsilon^{-1}$ follows immediately by noting that $\| \Upsilon^{-1} Y \|_{\ell^2}^2
= \sum_{n=1}^{\infty} \mu_n^{-2} |Y_n|^2 \le \mu_1^{-2} \| Y \|^2_{\ell^2}$ since $\mu_n \ge \mu_1$ for all $n \ge 1$.
The boundedness of $\tilde{\chi}$ follows by noting that
$$
(\tilde{\chi} \ast Y)_n = \sum_m \langle \psi_n , \chi \psi_m \rangle Y_m = \langle \psi_n , \chi {\mathcal{Y}} \rangle
$$
where ${\mathcal{Y}}(y,z) = \sum_m Y_m \psi_m(z,z)$.   Thus,    $(\tilde{\chi} \ast Y)_n$ is
the generalized Fourier coefficient of the function $\chi {\mathcal{Y}}$ and hence, by Parseval's
equality.
\begin{equation}
\sum_n |(\tilde{\chi} \ast Y)_n|^2 = \| \chi {\mathcal{Y}} \|^2_{L^2(\Omega)} \le \| \chi \|_{L^{\infty}}^2 \| {\mathcal{Y}} \|_{L^2(\Omega)}^2\ =  \| \chi \|_{L^{\infty}}^2 \| Y \|_{\ell^2}\   ,
\end{equation}
where the last step in this expression again used Parseval's equality
\end{Proof}


Intuitively, there are no linear terms in \eqref{eq:diag} in the equations for $\{a_k\}_{k=0}^N$ (except for the term $-A\sigma \langle \check{\chi}, \gamma\rangle$, which is decaying exponentially fast) or in the equation for $\sigma$. The equations for $\{b_k\}_{k=0}^N$ each contain a linear term of the form $-\Upsilon b_k$, where $\langle \Upsilon b_k, b_k \rangle \geq \mu_1 \|b_k\|^2$, with $\mu_1 > 0$. Hence, these variables should decay exponentially fast, and there is a spectral gap determined by $\mu_1$. Therefore, there should exist an invariant center-stable manifold of dimension $N+2$ of the form $\mathcal{M} = \{ b_k = h_k(a_0, a_1, \dots, a_N, \sigma): k = 0, 1, \dots N\}$. To see this, we note that $\gamma' = - \Upsilon \gamma$ and add this equation to \eqref{eq:diag} to obtain the autonomous system
\begin{eqnarray}
a_0' &=& 0 \nonumber \\
a_1' &=& -\frac{1}{2} \sigma^2 a_1 - A\sigma \langle \check{\chi}, \gamma \rangle_{\ell^2} \nonumber \\
a_k' &=& -\frac{k}{2}\sigma^2a_k - A\sigma^2  \langle \check{\chi}, b_{k-2} \rangle_{\ell^2} \qquad 2 \leq k \leq N \label{E:absg} \\
b_0' &=& - \Upsilon b_0 - A \tilde \chi \ast \gamma \nonumber \\
b_1' &=& -\left( \frac{1}{2} \sigma^2 + \Upsilon\right) b_1 - A\sigma \tilde \chi \ast \left[ b_0 - A a_0(\Upsilon^{-1}\check{\chi})\right] - D_{td} \sigma \gamma - \sigma A^2 \langle \check{\chi}, \gamma\rangle (\Upsilon^{-1}\check{\chi}) \nonumber \\
b_k' &=& -\left( \frac{k}{2}\sigma^2 + \Upsilon\right)  b_k - \sigma A \tilde \chi \ast \left[ b_{k-1} - Aa_{k-1}(\Upsilon^{-1}\check{\chi})\right]
- \sigma^2 D_{td} \left[ b_{k-2} - Aa_{k-2}(\Upsilon^{-1}\check{\chi})\right] \nonumber \\
&& \qquad   - \sigma^2A^2 \langle \check{\chi}, b_{k-2}\rangle (\Upsilon^{-1}\check{\chi}) \qquad \qquad 2 \leq k \leq N \nonumber \\
\sigma' &=& - \frac{1}{2} \sigma^3 \nonumber \\
\gamma' &=& -\Upsilon \gamma. \nonumber 
\end{eqnarray}
The linear part of this system (although no longer diagonal, due to the term $-A \tilde{\chi} \ast \gamma$ in the $b_0$ equation)
now makes the spectral separation clear. One could abstractly justify the existence of a center manifold of the form $(b_0, \dots b_N, \gamma) = H(a_0, \dots, a_N, \sigma)$. However, it turns out we can compute the function $H$ explicitly, and it has a rather simple form. Moreover, we can show directly that the center manifold is globally attracting. These results are collected in the following proposition.

\begin{Proposition} \label{prop:cm} For each $1 \leq k \leq N$, there exist functions $h_k = h_k(a_0, \dots, a_{k-1}, \sigma)$ of the form
\begin{equation}\label{E:defhk}
h_k(a_0, a_1, \ldots, a_{k-1}, \si) = \sum_{\ell=1}^{k} C^k_{k-\ell} a_{k-\ell} \si^\ell,
\end{equation}
where the $C^k_{k-\ell}$ are elements of $\ell^2$ for each $k$ and $\ell$, can be computed explicitly, are independent of $\nu$, and such that \eqref{E:absg} has an invariant center-stable manifold given by
\begin{equation}\label{E:defcm}
\mathcal{M}_N = \{ (b_0, \dots, b_N, \gamma) = (0, h_1(a_0, \sigma), \dots, h_N(a_0, \dots, a_{N-1}, \sigma), 0) \}.
\end{equation}
Moreover, there exist constants $C, \eta > 0$ that are independent of $\nu$ and such that all solutions to \eqref{E:absg} satisfy
\begin{equation}\label{E:decaytocm}
\|(b_0, \dots, b_N, \gamma)(T) - (0, h_1(a_0, \sigma), \dots, h_N(a_0, \dots, a_{N-1}, \sigma), 0)\|_{(\ell^2)^{N+2}} \leq C e^{-\eta T}, 
\end{equation}
where $(a_0, \dots, a_{N-1})$ and $\sigma$ are solutions of 
\begin{eqnarray*}
a_0' &=& 0 \nonumber \\
a_1' &=& -\frac{1}{2} \sigma^2 a_1 \nonumber \\
a_k' &=& -\frac{k}{2}\sigma^2a_k - A\sigma^2  \langle \check{\chi}, h_{k-2}(a_0, \dots, a_{k-3}, \sigma) \rangle_{\ell^2} \qquad 2 \leq k \leq N \\
\sigma' &=& - \frac{1}{2} \sigma^3. 
\end{eqnarray*}
Moreover, for all $k \geq 1$, 
\begin{equation}\label{E:cmdecay}
|a_k(\tau)| \leq C e^{-\eta \tau}, \qquad \tau = \log(1+T).
\end{equation}
\end{Proposition}

\begin{Remark} \label{rem:cmConst}
More precise statements of the convergence to the center manifold and decay within the center manifold are given in Lemmas \ref{lem:global} and \ref{lem:adecay}, respectively. Note that the exponential in $T$ convergence to the center manifold is equivalent to super-exponential in $\tau$ convergence, $e^{-\eta T} = e^{-\eta(e^\tau-1)}$, while the exponential in $\tau$ decay on the center manifold, implied by \eqref{E:cmdecay}, is equivalent to algebraic in $T$ decay, $e^{-\eta \tau} = (1+T)^{-\eta}$. Furthermore, the $\nu$- independence of the constants $C^k_{k-\ell}$ follows from the change of variable (\ref{E:nudep}).
\end{Remark}

\begin{Proof} 
The Proof will be divided into three steps: 1) Justifying \eqref{E:defhk}, the explicit formula for the center manifold; 2) Proving global convergence to the center manifold and justifying \eqref{E:decaytocm}; and 3) Justifying equation \eqref{E:cmdecay}, the decay rate within the center manifold.

{\bf Step 1: Explicit formula for the center manifold} To justify \eqref{E:defhk}, we will ultimately use induction, but we compute the first few terms directly since the equations in \eqref{E:absg} are different for $k = 0, 1$. First, notice that the set $(b_0, \gamma) = (0, 0)$ is invariant for \eqref{E:absg}. Next, we look for a function of the form
\[
h_1(a_0, \sigma) = C^1_0 a_0 \sigma, \qquad C^1_0 \in \ell^2,
\]
so that the set $(b_0, b_1, \gamma) = (0, h_1(a_0, \sigma), 0)$ is invariant. Computing $(b_0, b_1, \gamma)'$ in two different ways and equating the results, we find that we need
\[
-\frac{C^1_0}{2}a_0 \sigma^3 = -\frac{C^1_0}{2}a_0 \sigma^3 + \sigma a_0 \left[ - \Upsilon C^1_0 + A^2\tilde \chi \ast (\Upsilon^{-1}\check{\chi})\right].
\] 
Thus, we can take
\[
C^1_0 = A^2\Upsilon^{-1} \tilde \chi \ast (\Upsilon^{-1}\check{\chi}).
\]
Next, we look for a function of the form
\[
h_2(a_0, a_1, \sigma) = C^2_1 a_1 \sigma + C_0^2 a_0 \sigma^2
\]
so that the set $(b_0, b_1, b_2,  \gamma) = (0, h_1(a_0, \sigma), h_2(a_0, a_1,\sigma), 0)$ is invariant. As above, we find
\[
C_1^2 = A^2 \Upsilon^{-1} \tilde \chi \ast (\Upsilon^{-1} \check{\chi}), \qquad C_0^2 = D_{td} A \Upsilon^{-2} \check{\chi} - A^3 \Upsilon^{-1} [\tilde \chi \ast (\tilde \chi \ast(\Upsilon^{-1}\check{\chi}))].
\]
We now assume that \eqref{E:defhk} holds for $0 \leq k \leq n$ and prove this implies it is true for $k = n+1$ with $n \geq 2$. 
First, we compute
\begin{eqnarray}
b_{n+1}' &=& \frac{d}{dt} \sum_{\ell = 1}^{n+1} C_{n+1-\ell}^{n+1} a_{n+1-\ell} \sigma^\ell \nonumber \\
&=& \sum_{\ell = 1}^{n+1} \left( -\frac{1}{2} \right) \ell C_{n+1-\ell}^{n+1} a_{n+1-\ell} \sigma^{\ell+2} - \frac{1}{2} C_1^{n+1} a_1 \sigma^{n+2} - \sigma^2 C_2^{n+1} a_2 \sigma^{n-1} \nonumber \\
&& \qquad - \sum_{\ell =1}^{n-2}\frac{(n+1-\ell)}{2} C_{n+1-\ell}^{n+1}\sigma^{\ell+2} a_{n+1-\ell} - A \sum_{\ell =1}^{n-2} C_{n+1-\ell}^{n+1}\sigma^{\ell+2} \left \langle \check{\chi}, \sum_{j=1}^{n-\ell -1} C_{n-\ell-1-j}^{n-\ell-1} \sigma^j a_{n-\ell-1-j} \right\rangle. \nonumber \\
&=& -\frac{(n+1)}{2}  \sum_{\ell =1}^{n+1}C_{n+1-\ell}^{n+1}\sigma^{\ell+2} a_{n+1-\ell} - A \sum_{\ell =1}^{n-2} C_{n+1-\ell}^{n+1}\sigma^{\ell+2} \left \langle \check{\chi}, \sum_{j=1}^{n-\ell -1} C_{n-\ell-1-j}^{n-\ell-1} \sigma^j a_{n-\ell-1-j} \right\rangle.  \label{E:onehand}
\end{eqnarray}
Using \eqref{E:absg} and evaluating at $b_k = h_k$, $b_0 = \gamma = 0$, we also have
\begin{eqnarray}
b_{n+1}' &=& - \left[ \frac{(n+1)}{2}\sigma^2 + \Upsilon \right] \sum_{\ell = 1}^{n+1} C_{n+1-\ell}^{n+1} a_{n+1-\ell} \sigma^\ell - A \tilde \chi \ast \sum_{\ell = 1}^n C_{n-\ell}^n a_{n-\ell} \sigma^{\ell+1} \nonumber \\
&& \qquad - D_{td} \sum_{\ell = 1}^{n-1} C_{n-1-\ell}^{n-1} a_{n-1-\ell} \sigma^{\ell + 2} + \sigma A^2 a_n \tilde \chi \ast (\Upsilon^{-1} \check{\chi}) + D_{td} A \sigma^2 a_{n-1} (\Upsilon^{-1} \check{\chi}) \nonumber  \\
&& \qquad \qquad - A^2(\Upsilon^{-1} \check{\chi}) \left \langle \check{\chi}, \sum_{\ell = 1}^{n-1} C_{n-1-\ell}^{n-1} a_{n-1-\ell} \sigma^{\ell+2} \label{E:otherhand} \right\rangle.
\end{eqnarray}
We now equate the expressions on the right hand sides of equations \eqref{E:onehand}-\eqref{E:otherhand} to obtain
\begin{eqnarray*}
&& - A \sum_{\ell =1}^{n-2} C_{n+1-\ell}^{n+1}\sigma^{\ell+2} \left \langle \check{\chi}, \sum_{j=1}^{n-\ell -1} C_{n-\ell-1-j}^{n-\ell-1} \sigma^j a_{n-\ell-1-j} \right\rangle \\
&& \qquad = -\Upsilon  \sum_{\ell = 1}^{n+1} C_{n+1-\ell}^{n+1} a_{n+1-\ell} \sigma^\ell - A \tilde \chi \ast \sum_{\ell = 1}^n C_{n-\ell}^n a_{n-\ell} \sigma^{\ell+1} - D_{td} \sum_{\ell = 1}^{n-1} C_{n-1-\ell}^{n-1} a_{n-1-\ell} \sigma^{\ell + 2} \\
&& \qquad \qquad + \sigma A^2 a_n \tilde \chi \ast (\Upsilon^{-1} \check{\chi}) + D_{td} A \sigma^2 a_{n-1} (\Upsilon^{-1} \check{\chi})  
- A^2(\Upsilon^{-1} \check{\chi}) \left \langle \check{\chi}, \sum_{\ell = 1}^{n-1} C_{n-1-\ell}^{n-1} a_{n-1-\ell} \sigma^{\ell+2} \right\rangle.
\end{eqnarray*}
First, consider the resulting terms involving $a_n$. We need 
\[
0 = -\Upsilon C_n^{n+1} a_n \sigma + \sigma A^2 a_n \tilde \chi \ast (\Upsilon^{-1} \check{\chi}) \qquad \Rightarrow \qquad C_n^{n+1} = A^2 \Upsilon^{-1}\tilde \chi \ast (\Upsilon^{-1} \check{\chi}). 
\]
The terms involving $a_{n-1}$ imply
\[
0 = -\Upsilon C_{n-1}^{n+1}  - A \tilde \chi \ast C_{n-1}^n + D_{td} A  (\Upsilon^{-1}\check{\chi}) \qquad \Rightarrow \qquad C_{n-1}^{n+1} = D_{td} A \Upsilon^{-2} \check{\chi} - A \Upsilon^{-1}(\tilde \chi \ast C_{n-1}^n).
\]
The terms involving $a_{n-2}$ imply
\begin{eqnarray*}
0 = &&-\Upsilon C_{n-2}^{n+1}  - A \tilde \chi \ast C_{n-2}^n - D_{td} C_{n-2}^{n-1} - A(\Upsilon^{-1}\check{\chi}) \langle \check{\chi}, C_{n-2}^{n-1} \rangle \\
&& \qquad \Rightarrow \qquad C_{n-2}^{n+1} = \Upsilon^{-1} \left[  - A \tilde \chi \ast C_{n-2}^n - D_{td} C_{n-2}^{n-1} - A(\Upsilon^{-1}\check{\chi}) \langle \check{\chi}, C_{n-2}^{n-1} \rangle \right].
\end{eqnarray*}
Finally, for $3 \leq k \leq n$, the terms involving $a_{n-k}$ imply
\[
 - A \sum_{\ell =1}^{n-2} C_{n+1-\ell}^{n+1} \langle \check{\chi}, C_{n-k}^{n-\ell-1} \rangle =  -\Upsilon C_{n-k}^{n+1}  - A \tilde \chi \ast C_{n-k}^n - D_{td} C_{n-k}^{n-1} - A^2 (\Upsilon^{-1} \check{\chi}) \langle \check{\chi}, C_{n-k}^{n-1} \rangle,
\]
which gives
\[
C_{n-k}^{n+1} = \Upsilon^{-1} \left[ A \sum_{\ell =1}^{n-2} C_{n+1-\ell}^{n+1} \langle \check{\chi}, C_{n-k}^{n-\ell-1} \rangle- A \tilde \chi \ast C_{n-k}^n - D_{td} C_{n-k}^{n-1} - A^2 (\Upsilon^{-1} \check{\chi}) \langle \check{\chi}, C_{n-k}^{n-1} \rangle \right].
\]
All of the coefficients appearing in the sums on the RHS of this 
expression have been computed at previous stages of the iteration and 
hence we obtain $C^{n+1}_{n-k}$ in the form asserted in the Proposition.

{\bf Step 2: Proving global convergence to the center manifold and justifying \eqref{E:decaytocm}:} We'll show that the exact invariant manifolds previously constructed are globally attracting. First, note that we can solve \eqref{E:absg} explicitly to find
\begin{equation}\label{E:gammadecay}
\gamma_n(T) = \gamma_n(0) e^{- \mu_n T} \qquad \Rightarrow \qquad \|\gamma(T)\|_{\ell^2} \leq e^{- \mu_1 T}\|\gamma(0)\|_{\ell^2},
\end{equation}
\begin{equation}\label{E:b0decay}
b_0(T) = e^{- \Upsilon T} b_0(0) - \int_0^T e^{- \Upsilon(T-s)} A \tilde \chi \ast \gamma(s) \rmd s \qquad \Rightarrow \qquad \|b_0(T)\|_{\ell^2} \leq C(\|b_0(0)\|_{\ell^2}, \|\gamma(0)\|_{\ell^2}) (1+T) e^{- \mu_1 T},
\end{equation}
and
\begin{equation}\label{E:sigmabdd}
\sigma(T) = \frac{1}{\sqrt{T+1}} \qquad \Rightarrow \qquad |\sigma(T)| \leq 1. 
\end{equation}
Next, define 
\bea 
B_k = b_k - h_k(a_0, \dots, a_{k-1}, \sigma), \qquad k \geq 1,
\eea
where $h_k$ is defined in \eqref{E:defhk}.
\begin{Lemma} \label{lem:global}
There exists a $C > 0$, independent of $\nu$, such that for all $t > 0$, 
\begin{eqnarray*}
\|B_1 (T)\|_{\ell^2} &\leq& C (1+T)^{\frac{3}{2}}e^{-\mu_1 T} \\
\|B_k (T)\|_{\ell^2} &\leq& C (1+T)^{1 + \frac{k}{2}}e^{- \mu_1 T} \qquad  2 \leq k \leq N.
\end{eqnarray*}
\end{Lemma}
\begin{Proof}
For $k = 1$, we can compute $B_1'$ and solve the resulting equation explicitly to find
\begin{eqnarray*}
&&B_1(T) = e^{-\Upsilon T - \frac{1}{2}\log(T+1)} B_1(0) \\
&& \qquad - \int_0^T e^{- \Upsilon (T-s) - \frac{1}{2}(\log(T+1) - \log(s+1))} \left[ \frac{A}{\sqrt{1+s}} \tilde \chi \ast b_0(s) + \frac{D_{td}}{\sqrt{1+s}}\gamma(s) + \frac{A^2}{\sqrt{1+s}}(\Upsilon^{-1}\check{\chi}) \langle \check{\chi}, \gamma(s)\rangle \right] \rmd s.
\end{eqnarray*}
As a result, 
\begin{equation}\label{E:B1decay}
\|B_1(T)\|_{\ell^2} \leq C(\|B_1(0)\|_{\ell^2}, \|b_0(0)\|_{\ell^2}, \|\gamma(0)\|_{\ell^2}) (1+T)^{3/2} e^{-\mu_1 T}.
\end{equation}
Next, for $k \geq 2$, we have
\begin{eqnarray*}
B_k' &=& - \left( \frac{k^2}{2}\sigma^2 + \Upsilon\right) B_k - \sigma A \tilde \chi \ast B_{k-1} - \sigma^2 D_{td} B_{k-2} - \sigma^2 A^2 \langle \check{\chi}, B_{k-2} \rangle (\Upsilon^{-1} \check{\chi}),
\end{eqnarray*}
and so, assuming the result is true for $k-1$,
\begin{eqnarray*}
&&\|B_k(T)\|_{\ell^2} \leq e^{-\mu_1 T} \|B_k(0)\|_{\ell^2} \\
&& \qquad - C(\|B_{k-1}(0)\|_{\ell^2}, \|B_{k-2}(0)\|_{\ell^2}, \|\gamma(0)\|_{\ell^2})\int_0^T e^{-\mu_1 T} \left[ \frac{1}{\sqrt{1+s}} \frac{(1+s)^{1+\frac{k-1}{2}}}{1} + \frac{1}{(1+s)}\frac{(1+s)^{1+\frac{k-2}{2}}}{1}\right] \rmd s,
\end{eqnarray*}
which implies the result.
\end{Proof}

{\bf Step 3: Justifying equation \eqref{E:cmdecay}, the decay rate within the center manifold}

The goal of this section is to compute the decay rates of the $a_k$ by considering the system \eqref{eq:diag} reduced to its center manifold, which is given by
\begin{eqnarray*}
a_0' &=& 0 \nonumber \\
a_1' &=& -\frac{1}{2} \sigma^2 a_1 \nonumber \\
a_k' &=& -\frac{k}{2}\sigma^2a_k - A\sigma^2  \langle \check{\chi}, h_{k-2}(a_0, \dots, a_{k-3}, \sigma) \rangle_{\ell^2} \qquad 2 \leq k \leq N \\
\sigma' &=& - \frac{1}{2} \sigma^3. 
\end{eqnarray*}
Converting back to $\tau = \log(1+T)$, this becomes
\begin{eqnarray*}
\dot{a}_0 &=& 0 \nonumber \\
\dot{a}_1 &=& -\frac{1}{2}  a_1 \nonumber \\
\dot{a}_k &=& -\frac{k}{2} a_k - A  \langle \check{\chi}, h_{k-2}(a_0, \dots, a_{k-3}, e^{-\frac{\tau}{2}}) \rangle_{\ell^2} \qquad 2 \leq k \leq N. \end{eqnarray*}
Using the fact that $h_0 = 0$, we see immediately that
\begin{equation}\label{E:a012decay}
a_0(\tau) = a_0(0), \qquad a_1(\tau) = a_1(0)e^{-\frac{1}{2}\tau}, \qquad a_2(\tau) = a_2(0)e^{-\tau}.
\end{equation}

\begin{Lemma} \label{lem:adecay} There exists a $C > 0$, independent of $\nu$, such that if we write $k = 3j + n$ with $j \in \mathbb{N} \cup \{0\}$ and $n \in \{0, 1, 2\}$ then, for all $\tau \geq 0$, 
\begin{eqnarray*}
|a_k(\tau)| &\leq& C e^{ -\frac{(j+n)}{2} \tau}, \qquad 0 \leq k \leq N.
\end{eqnarray*}
\end{Lemma}

\begin{Proof}
Using the bound for $h_1$ in \eqref{E:defhk}, we find
\[
|a_3(\tau)| \leq e^{-\frac{3}{2}\tau} |a_3(0)| + C \int_0^t e^{-\frac{3}{2}(\tau-s)} |a_0(s)|e^{-\frac{s}{2}} \rmd s,
\]
which implies
\[
|a_3(\tau)| \leq Ce^{-\frac{1}{2}\tau}.
\]
A similar calculation shows
\[
|a_4(\tau)|  \leq C e^{-\tau}, \qquad |a_5(\tau)| \leq C e^{-\frac{3}{2}\tau}.
\]
Consider now general $k$, and assume the result holds for $a_m$ with $m \leq k-1$. Using \eqref{E:defhk}, we have
\[
|a_k(\tau)| \leq |a_k(0)|e^{-\frac{k}{2}\tau} + \int_0^\tau e^{-\frac{k}{2}(\tau-s)} \left( \sum_{\ell=1}^{k-2} C a_{k-2-\ell} e^{-\frac{\ell}{2}s}\right) \rmd s.
\]
Notice that
\[
\sum_{\ell=1}^{k-2} a_{k-2-\ell} e^{-\frac{\ell}{2}s} =  a_{k-3} e^{-\frac{1}{2}s} + a_{k-4} e^{-s} +  a_{k-5} e^{-\frac{3}{2}s} + \dots +  a_1 e^{-\frac{(k-3)}{2}s} + a_0 e^{-\frac{(k-2)}{2}s}.
\]
Thus, if $k = 3j + n$, we find
\[
\sum_{\ell=1}^{k-2}  a_{k-2-\ell} e^{-\frac{\ell}{2}s} \sim  e^{-\frac{(j+n)}{2}s} +  e^{-\frac{(j+n+2)}{2}s} + \dots +  e^{-\frac{(3j+n-2)}{2}s}.
\]
Thus, we find
\[
|a_k(\tau)| \leq |a_k(0)|e^{-\frac{k}{2}\tau} + \int_0^\tau e^{-\frac{k}{2}(\tau-s)} C e^{-\frac{(j+n)}{2}s} \rmd s \leq C e^{-\frac{(j+n)}{2}\tau}
\]
as claimed.
\end{Proof}

This concludes the proof of Proposition \ref{prop:cm}.
\end{Proof}


\subsection{Proof of Theorem \ref{thm:main1}(i)}\label{S:thm1i}

We now show how Theorem \ref{thm:main1}(i) follows from Proposition \ref{prop:cm}. Recall the definition of $u_\mathrm{app}$ in \eqref{E:defuapp}. The dynamics of $u_\mathrm{app}$ are governed by the behavior of $\{ \alpha_k\}_{k=0}^N$ and $\{\beta_k^n\}_{k=0}^N$, where $n = 1, 2, \dots$. Their dynamics are governed by \eqref{E:alphabeta}, which is a system of ODEs on $\mathbb{R}^N \times (\ell^2(\mathbb{R}))^N$. Proposition \ref{prop:cm} shows that, after converting to the variables $a_k, b_k^n$, this system has a finite-dimensional globally attracting center manifold given by \eqref{E:defcm}, and the rate of convergence to that center manifold is exponential in $T$, as given in \eqref{E:decaytocm}. Finally, recalling that $\alpha_k = a_k$, $\beta_k^n = b_k^n - (A\chi_n/\mu_n)\alpha_k$, and that the only term among $a_k$, $b_k^n$ that is not decaying in time is $a_0$, one obtains the leading behavior of \eqref{E:defuapp}. This justifies the statements in \ref{thm:main1}(i).


\section{Decay of the remainder via spectral decomposition and hypocoercivity}\label{S:remainder}

The goal of this section is to prove Theorem \ref{thm:main1}(ii), which states that the remainder terms decay rapidly. To that end, insert the expansion \eqref{E:deflowhigh} into \eqref{E:w0T3}-\eqref{E:vnT2} and project off the first $N+1$ eigenfunctions to obtain
\begin{eqnarray}
\partial_\tau w_0^s &=& \mathcal{L}_{td} w_0^s - D_{td}\left[ \alpha_{N-1} \varphi_{N+1}^{td} + \alpha_N \varphi_{N+2}^{td} + \partial_\xi^2 w_0^s \right] \nonumber \\
&& \qquad - A \sum_{m=1}^\infty \chi_m \left[ \beta_{N-1}^m \varphi_{N+1}^{td} + \beta_N^m \varphi_{N+2}^{td} + \partial_\xi^2 v_m^s \right] \label{E:rem} \\
\partial_\tau v_n^s &=& \mathcal{L}_{td} v_n^s - D_{td} \left[ \beta_{N-1}^n \varphi_{N+1}^{td} + \beta_N^n \varphi_{N+2}^{td} + \partial_\xi^2 v_n^s \right] - e^\frac{\tau}{2} A \sum_{m=1}^\infty \chi_{n,m} \left[ \beta_{N}^m \varphi_{N+1}^{td} + \partial_\xi v_m^s \right] \nonumber \\
&& \qquad - e^\tau[\mu_n v_n^s + A\chi_m w_0^s]. \nonumber
\end{eqnarray}
The operator $\mathcal{L}_{td}$, acting on $w_0^s$ and $v_n^s$, decays like $e^{-\frac{N+1}{2}\tau}$. In addition, the forcing terms in the above equation decay like $\alpha_k, \beta_k$ with $k \geq N-1$, which, due to Lemmas \ref{lem:global} - \ref{lem:adecay}, decay like $e^{-(j+n)\tau/2} \leq e^{-k\tau/6}$, for $k = 3j+n$. Therefore, we expect $w_0^s$ and $v_n^s$ to decay with the same rate as the forcing terms. 

To prove this, we will not work with the above system in the $(\xi, \tau)$ variables, but we will instead work in the Fourier space associated with the original $(X,T)$ variables. Using the fact that
\begin{equation}\label{E:cov}
u_0^s(X,T) = \frac{1}{\sqrt{T+1}} w_0^s\left( \frac{X}{\sqrt{T+1}}, \log(T+1)\right), \qquad u_n^s(X,T) = \frac{1}{(T+1)} \partial_\xi v_n^s\left( \frac{X}{\sqrt{T+1}}, \log(T+1)\right),
\end{equation}
we find
\begin{eqnarray*}
\partial_T u_0^s &=& \nu^2 \partial_X^2 u_0^s - A \sum_{m=1}^\infty \chi_m \partial_X u_m^s \\
&& \qquad - \frac{D_{td}}{(1+T)^{3/2}} \left[ \alpha_{N-1}(\log(T+1))\varphi_{N+1}^{td}\left (\frac{X}{\sqrt{T+1}}\right) + \alpha_{N}(\log(T+1))\varphi_{N+2}^{td}\left (\frac{X}{\sqrt{T+1}}\right)\right] \\
&& \qquad - \frac{A}{(1+T)^{3/2}} \sum_{m=1}^\infty \chi_m \left[ \beta_{N-1}^m(\log(T+1))\varphi_{N+1}^{td}\left (\frac{X}{\sqrt{T+1}}\right) + \beta_{N}^m(\log(T+1))\varphi_{N+2}^{td}\left (\frac{X}{\sqrt{T+1}}\right)\right] \\
\partial_T u_n^s &=& \nu^2 \partial_X^2 u_n^s - A \sum_{m=1}^\infty \chi_{n,m} \partial_X u_m^s - [\mu_n u_n^s + A \chi_n \partial_X u_0^s] \\
&& \qquad  - \frac{D_{td}}{(1+T)^2} \left[ \beta^n_{N-1}(\log(T+1))\varphi_{N+2}^{td}\left (\frac{X}{\sqrt{T+1}}\right) + \beta_{N}^n(\log(T+1))\varphi_{N+3}^{td}\left (\frac{X}{\sqrt{T+1}}\right)\right] \\
&& \qquad - \frac{A}{(1+T)^{3/2}} \sum_{m=1}^\infty \chi_{n,m} \beta_{N}^m(\log(T+1))\varphi_{N+2}^{td}\left (\frac{X}{\sqrt{T+1}}\right).\end{eqnarray*}
{\color{black} We now take the Fourier transform with respect to $x$, with the convention  
\[
\hat u(\kappa) = \int_{\mathbb{R}} e^{-\rmi \kappa x}u(x) \rmd x.
\]
Using} the notation
\begin{equation}
\hat U(\kappa, T) = \begin{pmatrix} \hat u_0^s(\kappa,T) \\ \{\hat u_n^s(\kappa,T)\}_{n=1}^\infty \end{pmatrix}, \qquad \check{\chi} = \{\chi_n\}_{n=1}^\infty, \qquad (\tilde \chi \ast f)_n = \sum_{m=1}^\infty \chi_{n,m} f_m, \label{E:defUhat}
\end{equation}
we find
\begin{equation}\label{E:remft}
\frac{d}{dT} \hat U = \mathcal{B}(\kappa) \hat U + \hat F(\kappa, T),
\end{equation}
where
\begin{equation}\label{E:defB}
\mathcal{B}(\kappa) = -\nu^2 \kappa^2 \begin{pmatrix} 1 & 0 \\ 0 & 1 \end{pmatrix} + \rmi \kappa A \begin{pmatrix} 0 & \check{\chi} \cdot \\ \check{\chi} & \tilde \chi \ast \end{pmatrix} - \begin{pmatrix} 0 & 0 \\ 0 & \Upsilon \end{pmatrix} =: \kappa^2 \mathcal{B}_2 + \kappa \mathcal{B}_1 + \mathcal{B}_0
\end{equation}
and
\begin{equation}\label{E:defF}
\hat F(\kappa, T) = \begin{pmatrix} \hat F_1(\kappa, T) \\ \hat F_2(\kappa, T) \end{pmatrix},
\end{equation}
with
\begin{eqnarray}
\hat F_1(\kappa,T) &=& - \frac{D_{td}\hat \Phi_0^{td}(\kappa,T)}{(1+T)^{3/2}} \left[ \alpha_{N-1}(T)(1+T)^{\frac{N+1}{2}}(-\rmi \kappa)^{N+1} + \alpha_{N}(T)(1+T)^{\frac{N+2}{2}}(-\rmi \kappa)^{N+2} \right] \nonumber  \\
&& \qquad - \frac{A\hat \Phi_0^{td}(\kappa,T)}{(1+T)^{3/2}} \sum_{m=1}^\infty \chi_m \left[ \beta_{N-1}^m(T)(1+T)^{\frac{N+1}{2}}(-\rmi \kappa)^{N+1}  + \beta_{N}^m(T)(1+T)^{\frac{N+2}{2}}(-\rmi \kappa)^{N+2} \right] \nonumber \\
\hat F_2(\kappa,T) &=& - \frac{D_{td}\hat \Phi_0^{td}(\kappa,T)}{(1+T)^2} \left[ \beta^n_{N-1}(T)(1+T)^{\frac{N+2}{2}}(-\rmi \kappa)^{N+2}+ \beta_{N}^n(T)(1+T)^{\frac{N+3}{2}}(-\rmi \kappa)^{N+3}\right] \nonumber \\
&& \qquad - \frac{A\hat \Phi_0^{td}(\kappa,T)}{(1+T)^{3/2}} \sum_{m=1}^\infty \chi_{n,m} \beta_{N}^m(T)(1+T)^{\frac{N+2}{2}}(-\rmi \kappa)^{N+2}. \label{E:defF2}
\end{eqnarray}
Note that we have written $\alpha_j(\log(T+1)) = \alpha_j (T)$ and $\beta_j^n(\log(T+1)) = \beta_j^n(T)$ for convenience, and a direct calculation shows that
\[
\hat \Phi_0^{td}(\kappa, T) = \sqrt{T+1}e^{-\nu_{td} \kappa^2 (T+1)}.
\]

The plan is to analyze the behavior of \eqref{E:remft} using Duhamel's formula,
\begin{equation}\label{E:remduhamel}
\hat U(\kappa,T) = e^{\mathcal{B}(\kappa)T} \hat U(\kappa, 0) + \int_0^T e^{\mathcal{B}(\kappa)(T-s)} F(\kappa, s) \rmd s, 
\end{equation}
and show that solutions decay like $T^{-\mathcal{N}(N)}$, where $\mathcal{N}$ can be made large by choosing $N$ large. The precise relationship between $\mathcal{N}$ and $N$ is given in the statement of Proposition \ref{prop:errordecay}. We will obtain this decay in the norm
\begin{equation}\label{E:Uhatnorm}
\|\hat U(\cdot, T) \|^2 = \int_\mathbb{R} \| \hat U(\kappa, T)\|_Y^2 \rmd \kappa= \int_\mathbb{R} |\hat u_0^\mathrm{s}(\kappa, T)|^2 \rmd \kappa + \int_\mathbb{R} \|\{\hat u_n^\mathrm{s}(\kappa, T)\}\|_{\ell^2}^2 \rmd \kappa. 
\end{equation}

\begin{Remark}
Recall that we expect decay of the remainder terms $w_0^s, v_n^s$ in $L^2(m)$, and the relationship between these variables and $u_0, u_n$ is given in \eqref{E:cov}. Suppose that two functions $g$ and $f$ are related via
\[
g(\xi, \tau) = (1+T)^\gamma f(X,T), \qquad \xi = \frac{X}{\sqrt{1+T}}, \qquad \tau = \log (1+T). 
\]
Then we have
\begin{eqnarray*}
\|g(\tau)\|_{L^2(m)}^2 &=& \int (1+\xi^2)^m |g(\xi, \tau)|^2 \rmd \xi = (1+T)^{2\gamma-1/2} \int [1 + X^2(1+T)^{-1}]^m|f(X,T)|^2 \rmd X \\
&\simeq& (1+T)^{2\gamma-1/2} \sum_{j=0}^m (1+T)^{-j} \int |X^j f(X,T)|^2 \rmd X \\
&=& (1+T)^{2\gamma-1/2} \sum_{j=0}^m (1+T)^{-j} \| \partial_\kappa^j \hat f(T) \|_{L^2}^2.
\end{eqnarray*}
The discussion at the beginning of this section suggests we can expect $w_0^\mathrm{s}(\xi, \tau)$ and $v_n^\mathrm{s}(\xi, \tau)$ to decay like
\[
\|w_0^\mathrm{s}(\tau)\|_{L^2(m)} + \| \|v^\mathrm{s}(\tau)\|_{\ell^2}\|_{L^2(m)} \sim e^{-\eta(N) \tau},
\]
where $\eta(N)$ grows with $N$. Therefore, one could estimate solutions to \eqref{E:remduhamel} in terms of the norm
\begin{equation}\label{E:triplenorm}
|||\hat U(T)||| = (1+T)^{1/2} \sum_{j=0}^m (1+T)^{-j} \| \partial_\kappa^j \hat u_0 (T) \|_{L^2}^2 + (1+T)^{3/2} \sum_{j=0}^m (1+T)^{-j} \| \partial_\kappa^j ( \kappa^{-1} \|\{\hat u_n (T)\}\|_{\ell^2}) \|_{L^2}^2.
\end{equation}
Although this is possible \cite{Chaudhary17}, the calculations are cumbersome. Therefore, we have chosen to carry out the estimates in terms of the much simpler norm \eqref{E:Uhatnorm}, which also seems quite natural.
\end{Remark}

The goal of this section will be to prove the following result.
\begin{Proposition}\label{prop:errordecay} For any $N \in \mathbb{N}$ and $\hat{U}(\kappa, 0)$ such that $\|\partial_\kappa^\ell \hat U(\cdot, 0)\| < \infty$ for all $0 \leq \ell \leq N+1$ and {\color{black} and $\partial_\kappa^\ell \hat U(0, 0) = 0$ for all $0 \leq \ell \leq N$} , the corresponding solution of \eqref{E:remduhamel} satisfies
\[
\| \hat U(\cdot, T) \| \leq C (1+T)^{-\frac{N}{6} - \frac{1}{12}}.
\]
for all $T \geq 0$, where $C$ is a constant that is independent of $\nu$ but depends on $\hat {U}(0)$ and its derivatives.
\end{Proposition}

\begin{Remark}
{\color{black} The assumption that $\partial_\kappa^\ell \hat U(0, 0) = 0$ for all $0 \leq \ell \leq N$ holds for initial data associated  with $\hat{U}(\kappa, \tau)$ defined in \eqref{E:defUhat}, due to equations \eqref{E:cov} and \eqref{E:deflowhigh} and the discussion following \eqref{E:defvarphi}. See also Lemma \ref{lem:TNzero}.}
\end{Remark}

\begin{Remark}\label{rem:thm1ii}
Note that the result claimed in Theorem \ref{thm:main1}(ii) follows from the above proposition. To see this, recall that $u_\mathrm{rem}$ is defined in \eqref{E:defurem}. Using equations \eqref{E:cov}, \eqref{E:Uhatnorm}, and Plancherel's Theorem , we have
\begin{eqnarray*}
\|u_\mathrm{rem}(T)\|_{L^2}^2 &\leq& C_\psi \left[ \|u_0^s(T)\|_{L^2}^2 +  \|\|\{u_n^s(T)\}\|_{\ell^2}^2\|_{L^2}^2 \right] \\
&=& C_\psi \left[ \|\hat u_0^s(T)\|_{L^2}^2 +  \|\|\{\hat u_n^s(T)\}\|_{\ell^2}^2\|_{L^2}^2 \right] \\
&=& C_\psi \|\hat U(T)\|^2 \leq C (1+T)^{-\frac{N}{3} - \frac{1}{6}},
\end{eqnarray*}
where $C_\psi$ is a constant that depends on the $L^2$ norms of the cross-sectional eigenfunctions $\psi_n$.
Note that the requirement that $\|\partial_\kappa^\ell \hat U(\cdot, 0)\| < \infty$ for all $0 \leq \ell \leq N+1$ in the above proposition holds as long as the initial data for \eqref{E:orig-model} lies in the algebraically weighted function space: $u(\cdot,0) \in L^2(N+1)\times L^2(\Omega)$. This is because $\partial_\kappa^\ell \hat f \in L^2$ if and only if $X^\ell f \in L^2$, which means $f \in L^2(\ell)$.
\end{Remark}

We now state a brief result on the decay of the forcing terms in \eqref{E:remduhamel}.

\begin{Lemma}\label{lem:Fbound}
There exists a constant $C > 0$, independent of $\nu$, such that, for all $T > 0$, $\kappa \in \mathbb{R}$
\begin{eqnarray*}
|\hat F_1(\kappa,T)| &\leq& C (1+T)^{\frac{N-1}{2} - \frac{1}{2}(j+n)} |\kappa|^{N+1} e^{-\nu_{td}\kappa^2(1+T)} [1 + |\kappa| (1+T)^{1/2}] \\
\|\hat F_2(\kappa,T)\|_{\ell^2} &\leq& C (1+T)^{\frac{N-1}{2} - \frac{1}{2}(j+n)} |\kappa|^{N+2} e^{-\nu_{td}\kappa^2(1+T)} [1 + | \kappa|(1+T)^{1/2}+ (1+T)^{1/2}]
\end{eqnarray*}
where $n, j$ are defined so that $N-1 = 3j+n$, with $n \in \{0, 1, 2\}$.
\end{Lemma}

\begin{Proof}
This is a direct consequence of the definition of $\hat F$ in \eqref{E:defF2}, of Lemmas \ref{lem:global} - \ref{lem:adecay}, and of \eqref{eq:diagchange}.
\end{Proof}

In order to combine Lemma \ref{lem:Fbound} with equation \eqref{E:remduhamel} and prove a decay result for the remainder terms, we will need good control of the semigroup generated by $\mathcal{B}(\kappa)$. To obtain this, we will first obtain estimates on the spectrum of $\mathcal{B}(\kappa)$. We will then use these spectral estimates to obtain decay estimates on the semigroup for three different regions: 1) small wavenumber $0 \leq |\kappa| \leq \kappa_0$; 2) intermediate wavenumber $\kappa_0  \leq |\kappa| \leq \kappa_1 \nu^{-1}$; and 3) large wavenumber $\kappa_1 \nu^{-1} \leq |\kappa|$, where $\kappa_0$ and $\kappa_{1}$ are positive constants that are independent of $\nu$.


\subsection{Spectral decomposition}\label{S:spectral}

First, we state a lemma on the spectrum of $\mathcal{B}_{0,1,2}$.
\begin{Lemma} On the space $Y = \mathbb{C} \times \ell^2(\mathbb{C})$ the following hold.
\begin{enumerate}
\item The operator $\mathcal{B}_0$ has only point spectrum, and it is given by $\sigma(\mathcal{B}_0) = \{0\} \cup \{-\mu_n\}_{n=1}^\infty$.
\item The operators $\mathcal{B}_1$ and $\mathcal{B}_2$ are bounded. 
\end{enumerate}
\end{Lemma}
\begin{Proof}
\begin{enumerate}
\item This follows from the fact that $\mathcal{B}_0$ is diagonal and the only accumulation point of its entries is $\infty$. 
\item This is trivially true for $\mathcal{B}_2$ because it is a scalar multiple of the identity, and for $\mathcal{B}_1$ it follows from the fact that $\{\psi_n\}_{n=0}^\infty$ forms an orthonormal basis for $L^2(\Omega)$ and Parseval's identity. 
\end{enumerate}
\end{Proof}

Next, we analyze the spectrum of $\mathcal{B}(\kappa)$ for any fixed $\kappa \in \mathbb{R}$. 
\begin{Lemma} \label{lem:point}
Fix any $\kappa \in \mathbb{R}$. The spectrum of $\mathcal{B}(\kappa)$ consists only of point spectrum.
\end{Lemma}
\begin{Proof}
We will show that, for fixed $\kappa$, $\mathcal{B}(\kappa) = \mathcal{B}_0 + \kappa(\mathcal{B}_1 + \kappa \mathcal{B}_2)$ is a relatively compact perturbation of $\mathcal{B}_0$. The result will then follow from Weyl's theorem \cite[XIII.4, Corollary 2]{ReedSimon78}. We must show
\[
\kappa (\mathcal{B}_1 + \kappa \mathcal{B}_2) (\mathcal{B}_0 + i)^{-1}
\]
is a compact operator on $\mathbb{C} \times \ell^2(\mathbb{C})$. By Parseval's identity, this is equivalent to showing that
\[
\kappa(\rmi A \chi(y,z) - \nu^2 \kappa) ( \Delta + i)^{-1}
\]
is a compact operator on $L^2(\Omega)$. We let $\{\hat{u}_n(y,z)\} \subset L^2(\Omega)$ be a bounded sequence: $||\hat{u}_n(y,z)||_{L^2(\Omega)} \leq C$ for all $n \in \mathbb{N}$. Then, since $i$ is in the resolvent set of $ \Delta$ and $( \Delta + i)^{-1}:L^2(\Omega) \rightarrow H^1(\Omega)$ is bounded, it follows that $\{(\Delta + i)^{-1} \hat{u}_n \}$ is a bounded sequence in $H^1(\Omega)$. Therefore 
\beas  
\{\kappa(\rmi A \chi(y,z) - \nu^2 \kappa)( \Delta + i)^{-1} \hat{u}_n \}
\eeas  
is also a bounded sequence in $H^1(\Omega)$. Since $H^1(\Omega)$ is compactly embedded in $L^2(\Omega)$, this sequence has an $L^2(\Omega)$ convergent subsequence. Therefore $\kappa( \rmi A \chi(y,z) - \nu^2 \kappa) ( \Delta + i)^{-1}$ is compact. \end{Proof}


\subsubsection{Low wavenumber estimates using the leading eigenvalue}

We next prove a result on the spectrum of $\mathcal{B}(\kappa)$ for $|\kappa|$ sufficiently small. In particular, we show in this case that the eigenvalues of $\mathcal{B}(\kappa)$ split into two parts: an eigenvalue $\lambda_0(\kappa)$ near $0$, and eigenvalues $\lambda(\kappa)$ satisfying $Re(\lambda(\kappa))\leq - \mu_1/2$. Therefore, we expect $\lambda_0(\kappa)$ to dominate the long-time behavior, and we will therefore be able to use it to obtain estimates on the low-wavenumber part of our solution. In addition, we will show that this leading eigenvalue $\lambda_0(\kappa)$ is approximately $-\nu_{td} \kappa^2$, so the long-time behavior will correspond with Taylor dispersion. 

{\color{black} We note that, at various points in the following proofs we will need to fix a constant $\kappa_0$ that is sufficiently small and consider only $\kappa$ such that $|\kappa| \leq \kappa_0$. The value of $\kappa_0$ will always be independent of $\nu$ and will only be adjusted a finite number of times. 
}

\begin{Proposition}\label{prop:spec-decomp}

{\color{black} There exists a sufficiently small constant $\kappa_0$ that is independent of $\nu$ and such that the following holds.} {\color{black} Fix any $\kappa \in \mathbb{R}$ such that $|\kappa| \leq \kappa_0$, and let $0< \nu < 1$.
\begin{enumerate}
\item The (point) spectrum of $\mathcal{B}(\kappa)$ can be divided into two disjoint sets, $\sigma(\mathcal{B}(\kappa)) = \{ \lambda_0(\kappa)\} \cup \Sigma(\kappa)$, where $|\lambda_0(\kappa)+\nu \kappa^2| \leq \sqrt{2}\mu_1/2$ and, for any eigenvalue $\lambda(\kappa) \in \Sigma(\kappa)$, we have $\mathrm{Re}(\lambda(\kappa)) \leq - \mu_1/2$.
\item The leading eigenvalue satisfies $\lambda_0(\kappa) = -\nu_{td} \kappa^2 + \Lambda_0(\kappa)$, where $\Lambda_0(\kappa) = \rmi r \kappa^3 + \mathcal{O}(\kappa^4)$ is smooth, and independent of $\nu$. Here $r = r(\chi, \{\mu_n\}_{n=1}^\infty) \in \mathbb{R}$ is given in equation \eqref{E:defr}. 
\end{enumerate}
}
\end{Proposition}

{\color{black}
The main idea behind this Proposition is the following: recall that $\mathcal{B}(\kappa) = \mathcal{B}_0 + \kappa \mathcal{B}_1 + \kappa^2 \mathcal{B}_2$. If $|\kappa|$ is small, then $\mathcal{B}(\kappa)$ is just a small perturbation of $\mathcal{B}_0$, which has spectrum $\{0\} \cup \{-\mu_n\}_{n=1}^\infty$ and the separation claimed in (i). Furthermore, we will see that $\mathcal{B}_1$ is antisymmetric, hence the real part of the spectrum of $\mathcal{B}(\kappa)$ is actually an $\mathcal{O}(\kappa^2)$ perturbation of that of $\mathcal{B}_0$. The $\nu$-dependence of the spectrum stated in the proposition can be obtained from the following decomposition: recall that $\mathcal{B}_2 = - \nu^2 I$. Letting $\mathcal{C}(\kappa) = B_0 + \kappa B_1$, we have that $\mathcal{B}(\kappa) = \mathcal{C}(\kappa) - \nu^2 \kappa^2 I$. That is, the operators $\mathcal{B}(\kappa)$ and $\mathcal{C}(\kappa)$ differ by a scalar multiple of the identity, and, since $\mathcal{C}(\kappa)$ is independent of $\nu$, all of the $\nu$-dependence of $\mathcal{B}(\kappa)$ is contained in this scalar. Therefore we immediately have the following lemma:
\begin{Lemma}\label{prop:spec-nu}
Fix any $\kappa \in \mathbb{R}$, let $\nu > 0$, and let $\mathcal{B}(\kappa)$ and $\mathcal{C}(\kappa)$ be defined as above. The following are true:
\begin{enumerate}
\item The semigroups of $\mathcal{B}(\kappa)$ and $\mathcal{C}(\kappa)$ are related by $e^{\mathcal{B}(\kappa)T} = e^{-\nu^2 \kappa^2 T}e^{\mathcal{C}(\kappa)T}$.
\item The eigenvalues $\lambda(\kappa)$ of $\mathcal{B}(\kappa)$ and $\Gamma(\kappa)$ of $\mathcal{C}(\kappa)$ are in one-to-one correspondence with one another via $\lambda(\kappa) = \Gamma(\kappa) - \nu^2 \kappa^2$, and corresponding eigenvalues have the same projection operators $P(\kappa)$.
\end{enumerate}
\end{Lemma}
\begin{Remark} \label{rmk:spec-nu}
Since the operator $\mathcal{C}(\kappa)$ is independent of $\nu$, the above lemma tells us exactly what the $\nu$-dependence is in the semigroup $e^{\mathcal{B}(\kappa)T}$, and it tells us exactly what the $\nu$-dependence is in the eigenvalues $\lambda(\kappa)$ in terms of the ($\nu$-independent) eigenvalues $\Gamma(\kappa)$ of $\mathcal{C}(\kappa)$. Furthermore, since the projections $P(\kappa)$ of corresponding eigenvalues are the same, and $\mathcal{C}(\kappa)$ is independent of $\nu$, these projections can be taken to be independent of $\nu$. This relationship between the $\nu$-dependence and the structure of the system is a direct consequence of the change of variables (\ref{E:nudep}).
\end{Remark}

Note that, because $\mathcal{B}(\kappa)$ generates an analytic semigroup, the following Corollary follows immediately from Proposition \ref{prop:spec-decomp}(i). 

\begin{Corollary} \label{cor:explow}
{\color{black} There exists a sufficiently small constant $\kappa_0$ that is independent of $\nu$ and such that the following holds.} Fix any $\kappa \in \mathbb{R}$ such that $|\kappa| \leq \kappa_0$, and let $0< \nu < 1$. Let $Q_0(\kappa)$ be the projection complementary to the eigenspace of the eigenvalue $\lambda_0(\kappa)$ of $\mathcal{B}(\kappa)$. Then, for all $W \in \mathbb{C} \times \ell^2(\mathbb{C}) = Y$ and $T > 0$, we have 
\[
\| e^{\mathcal{B}(\kappa)T} Q_0(\kappa) W \|_Y \leq C e^{-\frac{ \mu_1}{2} T}\|W\|_Y,
\] 
for some constant $C > 0$ which is independent of $\nu$.
\end{Corollary}

Before proving Proposition \ref{prop:spec-decomp}, we will need to prove the following Lemma. 

\begin{Lemma} \label{lem:specBoundSmallk}
{\color{black} There exists a sufficiently small constant $\kappa_0$ that is independent of $\nu$ and such that the following holds.}  Let $\lambda(\kappa)$ be an eigenvalue of $\mathcal{B}(\kappa)$. Then
\begin{enumerate}
\item $Re(\lambda(\kappa)) \leq - \nu^2 \kappa^2$.
\item If $|\kappa| \leq \kappa_0$, then $|Im \lambda(\kappa)| < \mu_1/2$.
\end{enumerate}
\end{Lemma}

\begin{Proof} This lemma follows by splitting $\mathcal{B}(\kappa)$ into its real and imaginary parts.
Recall from \eqref{E:defB} that $\mathcal{B}(\kappa) = \mathcal{B}_0 + \kappa \mathcal{B}_1 + \kappa^2 \mathcal{B}_2$ with
\[
\mathcal{B}_0 = \bpm 0 && 0 \\ 0 &&  -\Upsilon \epm, \qquad 
\mathcal{B}_1 = A i \bpm 0 && \check{\chi} \cdot  \\ \check{\chi} && \tilde \chi \ast \epm, \qquad 
\mathcal{B}_2 = -\nu^2 \bpm 1 && 0 \\ 0 && I \epm. 
\]
Note that $\mathcal{B}_0$ and $\mathcal{B}_2$ are diagonal and hence
\[
\mathcal{S}(\kappa):=\mathcal{B}_0 + \kappa^2 \mathcal{B}_2 
\]
is symmetric. Also note that 
\[ 
\mathcal{A}(\kappa):=\kappa \mathcal{B}_1
\] 
is anti-symmetric, which follows from a straightforward computation using Parseval's identity.  
Let $V=\{V_n\}_{n=0}^{\infty} \in \mathbb{C} \times \ell^2(\mathbb{C})$ and let $v(y,z) = V_0 + \sum_{n=1}^{\infty} V_n \psi_n(y,z)$.  Then
\begin{eqnarray*}
\langle \mathcal{B}_1 V, V \rangle_{\mathbb{C} \times \ell^2(\mathbb{C})} &=& \langle A i \chi(y,z) v(y,z), v(y,z) \rangle_{L^2(\Omega)}  \\
&=&  - \langle v(y,z), A i \chi(y,z) v(y,z) \rangle_{L^2(\Omega)} \\
& = & - \langle V, \mathcal{B}_1 V \rangle_{\mathbb{C} \times \ell^2(\mathbb{C})}.
\end{eqnarray*}

Using this splitting into symmetric and antisymmetric parts, if  $\lambda(\kappa)$ is an eigenvalue of $\mathcal{B}(\kappa)$ with eigenvector $V(\kappa)$ normalized so that $\|V(\kappa)\|_Y = 1$, one can immediately write (see \cite{Gallay:2005}, p. 124, for example)
\begin{equation}
Re(\lambda(\kappa)) = \langle {V(\kappa) } , \mathcal{S}(\kappa) V(\kappa) \rangle_{ \mathbb{C} \times \ell^2(\mathbb{C})}
\end{equation}
Since $\mathcal{S}(\kappa)$ is symmetric and $V(\kappa)$ is normalized, the variational
characterization of the eigenvalues of symmetric
operators insures that this inner product is bounded by the right-most point in the spectrum
$\mathcal{S}(\kappa)$ which -$\nu^2 \kappa^2$.  This
proves the first part of Lemma \ref{lem:specBoundSmallk}. 

For the second part of this lemma, we use an argument similar to that used in the proof of the first part to control the imaginary part of $\lambda(\kappa)$.  Writing 
$\lambda(\kappa) =  \langle {V(\kappa) } , \mathcal{B}(\kappa) V(\kappa) \rangle_{ \mathbb{C} \times \ell^2(\mathbb{C})}$
and splitting $\mathcal{B}$ into its symmetric and anti-symmetric parts
 yields an expression for $Im (\lambda(\kappa))$: 
\beas 
Im (\lambda(\kappa)) = -\frac{1}{area(\Omega)} \int_{\Omega} v(\kappa,y,z) \overline{\kappa i A \chi(y,z) v(\kappa,y,z)} dy dz, 
\eeas
where $v(\kappa,y,z)$ is the unit eigenvector for $\lambda(\kappa)$, and we have used Parseval's identity. Continuing, we get
\beas 
Im (\lambda(\kappa)) = i \kappa A \frac{1}{area(\Omega)} \int_{\Omega} v(\kappa,y,z) \chi(y,z) \overline{ v(\kappa,y,z)} dy dz, 
\eeas
so that 
\beas
|Im \lambda(\kappa) | \leq |\kappa| A \| \chi \|_{L^{\infty}_\Omega} < \frac{\mu_1}{2},
\eeas 
{\color{black} as long as $\kappa_0 < \mu_1/(2  A \| \chi \|_{L^{\infty}_\Omega} )$.}
This completes the proof of the Lemma.
\end{Proof}

We now prove Proposition \ref{prop:spec-decomp}
}
\begin{Proof} 
{\color{black}First, we prove item (i).  To establish this separation for
$\mathcal{B}(\kappa)$, we note first that it suffices to establish it for $\mathcal{C}(\kappa)$
since $\kappa^2 \mathcal{B}(\kappa)$ simply shifts the entire spectrum by an amount $\kappa^2$.
Let $\Gamma_*$ be the boundary of the rectangle $\{z = x + i y : |x|, |y| \leq \mu_1/2\}$. The
$\Gamma_* $ separates the spectrum of $\mathcal{B}_0$, and for $z \in \Gamma_*$ we have
\beas 
||(\mathcal{B}_0-z)^{-1}|| \leq \frac{2}{ \mu_1}\ ,
\eeas 
since $\mathcal{B}_0$ is diagonal and $ \mu_1/2$ is the distance from $\Gamma_*$ to $\si(\mathcal{B}_0)$ \cite[Chapter V, \S3.5, (3.16)]{Kato95}. 

Next note that the norm of $\kappa \mathcal{B}_1$, considered as an operator on
$\mathbb{C} \times \ell^2(\mathbb{C})$, can be bounded with the aid of Parseval's identity.  Let
$V \in \mathbb{C} \times \ell^2(\mathbb{C})$ and set $v(y,z,\kappa) = V_0(\kappa) + \sum_{n=1}^{\infty} V_n(\kappa) \psi_n(y,z)$.  Then
\begin{eqnarray}
&& \| \kappa \mathcal{B}_1 V \|^2_{\mathbb{C} \times \ell^2(\mathbb{C})} = 
|\kappa|^2 \langle A i \chi v,  A i \chi v \rangle \le |\kappa |^2 A^2 \| \chi \|_{L^\infty(\Omega)}^2 \| v \|_{L^2(\Omega)}^2 \\ \nonumber 
&& \qquad \qquad  = |\kappa |^2 A^2 \| \chi \|_{L^\infty(\Omega)}^2 \| V \|_{\mathbb{C} \times \ell^2(\mathbb{C})}^2
\end{eqnarray}
Thus, {\color{black} if $\kappa_0 <  \mu_1/(2  A \| \chi \|_{L^{\infty}_\Omega} )$,} $\| \kappa \mathcal{B}_1 \| \le \mu_1/2$.  

This in turn implies that for any $z \in \Gamma_* $,
\begin{equation}
\| (\mathcal{C}(\kappa) - z)^{-1} \| = \| (\mathcal{B}_0 + \kappa \mathcal{B}_1 -z)^{-1}\|
= \| ({ \bf{1}}  + \kappa (\mathcal{B}_0 -z)^{-1} \mathcal{B}_1 )^{-1} (\mathcal{B}_0 -z)^{-1} \|\  .
\end{equation}
By the estimate of the norm of $\mathcal{B}_1$ and the assumption that $|\kappa| < \kappa_0$,
we see that 
\begin{equation}
\|  \kappa (\mathcal{B}_0 -z)^{-1} \mathcal{B}_1 \| < 1\ ,
\end{equation}
so that $({ \bf{1}}  + \kappa (\mathcal{B}_0 -z)^{-1} \mathcal{B}_1 )^{-1}$ is bounded and
hence that $\Gamma_* $ is contained in the resolvent set of $\mathcal{C}(\kappa)$ for all
$| \kappa | \le \kappa_0$. Since the eigenvalues of $\mathcal{C}(\kappa)$ vary continuously
with $\kappa$, this means that there is one eigenvalue, $\Gamma_0(\kappa)$, of 
$\mathcal{C}(\kappa)$ inside $\Gamma_*$ for all $| \kappa | \le \kappa_0$, and hence that 
$|\Gamma_0(\kappa)| \le \sqrt{2}  \mu_1/2$. As we observed above, the corresponding eigenvalue of $\mathcal{B}(\kappa)$ is $\lambda_0(\kappa) = \Gamma_0(\kappa) - \nu \kappa^2$, and hence the first part of point $(i)$ in the Proposition follows.

Now suppose that $\Gamma(\kappa)$ is an eigenvalue not contained in $\Gamma_*$ (and hence, by the relationship between the spectra of ${\mathcal{B}}(\kappa)$ and ${\mathcal{C}}(\kappa)$ it corresponds to an eigenvalue $\lambda(\kappa)
\in \Sigma(\kappa)$).  Then either
\begin{itemize}
\item[(a)] $\Re(\Gamma(\kappa)) \le - \mu_1/2$\ , or
\item[(b)] $-\mu_1/2 \le \Re(\Gamma(\kappa)) \le 0$, and $| \Im(\Gamma(\kappa))| > \mu_1/2$,
\end{itemize}
because Lemma \ref{lem:specBoundSmallk} implies that none of the eigenvalues of ${\mathcal{C}}(\kappa)$ can have positive real part. If case (b) held, then there would be a corresponding eigenvalue $\lambda(\kappa)$ of ${\mathcal{B}}(\kappa)$ with
$| Im(\lambda(\kappa)) | > \mu_1/2$, and this would violate Lemma \ref{lem:specBoundSmallk} (ii).  Hence
case (a) applies and this in turn implies the bound in Proposition \ref{prop:spec-decomp} (i).
}

Next, we prove item (ii) in Proposition \ref{prop:spec-decomp}. Note that, because $\lambda_0(\kappa)$ is a perturbation of the simple eigenvalue $0$ of $\mathcal{B}_0$, both $\lambda_0(\kappa)$ and its spectral projection $P_0(\kappa)$ perturb smoothly in $\kappa$ \cite{Kato95}. However, due to Lemma \ref{prop:spec-nu}, we can instead estimate the leading ($\nu$-independent) eigenvalue $\Gamma_0(\kappa)$ of $\mathcal{C}(\kappa)$, which is still a perturbation of the simple eigenvalue $0$ of $\mathcal{B}_0$. We expand this eigenvalue
\bea \label{eq:evalexp}
\Gamma_0(\kappa) = \Gamma_0 + \Gamma_1 \kappa + \Gamma_2 \kappa^2 + \mathcal{O}(\kappa^3)
\eea 
and its corresponding eigenvector
\bea \label{eq:evectexp}
\hat{V}(\kappa) = \hat{V}_0 + \hat{V}_1 \kappa + \hat{V}_2 \kappa^2 + \mathcal{O}(\kappa^3),
\eea 
where
\[
\hat{V}(\kappa)= \begin{pmatrix} \hat{u}_0(k) \\ \hat{U}(\kappa) \end{pmatrix}, \qquad \hat{V}_j= \begin{pmatrix} \hat{u}_0^j \\ \hat{U}^j\end{pmatrix}.
\]
Now the eigenvalue problem reads 
\bea \label{eq:evalprob}
\mathcal{C}(\kappa) \hat{V}(\kappa) = \Gamma(\kappa) \hat{V}(\kappa), 
\eea  
Plugging \eqref{eq:evalexp} and \eqref{eq:evectexp} into \eqref{eq:evalprob}, we find
\[
\mathcal{B}_0 \hat{V}_0 =  0 \cdot \hat{V}_0, \qquad \Rightarrow \qquad \Gamma_0 = 0, \qquad \hat{V}_0 = \begin{pmatrix} 1 \\ 0 \end{pmatrix},
\]
Next, we find
\begin{align} \label{eq:perturbeqns}
\begin{split} 
\mathcal{B}_1 \hat{V}_0 + \mathcal{B}_0 \hat{V}_1 & =  \Gamma_1 \hat{V}_0 \\
\mathcal{B}_1 \hat{V}_1 + \mathcal{B}_0 \hat{V}_2 & =  \Gamma_2 \hat{V}_0 + \Gamma_1 \hat{V}_1 \\
\mathcal{B}_1 \hat{V}_2 + \mathcal{B}_0 \hat{V}_3 & =  \Gamma_3 \hat{V}_0 + \Gamma_2 \hat{V}_1 + \Gamma_1 \hat V_2
\end{split}
\end{align}
and so on. Solving the first equation, we find
\[
\Gamma_1 = 0, \qquad \hat V_1 = \begin{pmatrix} c_1 \\  \rmi A\Upsilon^{-1} \check{\chi} \end{pmatrix},
\]
where the scalar constant $c_1$ is undetermined but can be fixed by normalizing the eigenvectors. At $\mathcal{O}(\kappa^2)$, we similarly find
\[
\Gamma_2 = - D_{td}, \qquad \hat V_2 = \begin{pmatrix} c_2 \\ \rmi Ac_1 \Upsilon^{-1}\check{\chi} - A^2 \Upsilon^{-1}[\tilde \chi\ast(\Upsilon^{-1}\check{\chi})] \end{pmatrix}.
\]
Finally, at $\mathcal{O}(\kappa^3)$, the first component in the equation implies
\begin{eqnarray*}
\Gamma_3 &=& c_1(D_{td}) + \rmi A \check{\chi} \cdot \left[ \rmi Ac_1 \Upsilon^{-1}\check{\chi} - A^2[\tilde \chi \ast(\Upsilon^{-1}\check{\chi})] \right] \\
&=& - \rmi A^3 \check{\chi} \cdot [\tilde \chi \ast(\Upsilon^{-1}\check{\chi})]. 
\end{eqnarray*}
In particular, $\Gamma_3$ is purely imaginary, and therefore 
\[
\Gamma_0(k) = -D_{td} \kappa^2 + i r \kappa^3 + \mathcal{O}(\kappa^4),
\]
where
\begin{equation}\label{E:defr}
r = -A^3 \check{\chi} \cdot [\tilde \chi \ast(\Upsilon^{-1}\check{\chi})].
\end{equation}
Finally, using Lemma \ref{prop:spec-nu}, we have
\beas
\lambda_0(\kappa) & = & -(\nu^2 + D_{td})\kappa^2 + \Lambda_0(\kappa) \\
 & = & -\nu_{td} \kappa^2 + \Lambda_0(\kappa),
\eeas
where $\Lambda_0(\kappa) = i r \kappa^3 + \mathcal{O}(\kappa^4)$ is independent of $\nu$.
This completes the proof of item (ii), and of Proposition \ref{prop:spec-decomp}. 
\end{Proof}


\subsubsection{High wavenumber estimates using standard diffusive estimates}

Next, we consider the behavior of the spectrum of $\mathcal{B}(\kappa)$ for large $|\kappa|$. 
\begin{Corollary} \label{cor:exphigh} Given any fixed constant $\kappa_1$, for all $|\kappa| \geq \frac{\kappa_1}{\nu}$ we have
\[
\|e^{\mathcal{B}(\kappa)T} W \|_Y \leq C e^{- \kappa_1^2 T}\|W\|_Y.
\]
\end{Corollary}
\begin{Proof} This follows immediately from Lemma \ref{lem:specBoundSmallk}, using the fact that $\mathcal{B}(\kappa)$ generates an analytic semigroup. 
\end{Proof}


\subsubsection{Intermediate wavenumber estimates via hypocoercivity}

In this subsection, we prove the following Lemma. 

\begin{Proposition} \label{lem:intermediate-wn}
{\color{black} There exists a constant $\kappa_0$ sufficiently small and independent of $\nu$ so that the following holds.}
There exist positive constants $\kappa_1$ and $\delta \in (0,\frac{1}{4})$ such that for all $ \frac{\kappa_1}{\nu} \geq |\kappa| \geq \kappa_0(1-\delta)$ and $T>0$, we have
\[
\| e^{\mathcal{B}(\kappa) T} W\|_Y \leq C e^{- M T}\|W\|_Y,
\]
where $M$ and $C$ are {\color{black}positive constants that are} independent of $\nu$ {\color{black} and $\kappa$}.
\end{Proposition}

\begin{Remark}
This result does not appear to be obvious. A naive estimate, such as that in the proof of Corollary \ref{cor:exphigh}, would only give
\[
\|e^{\mathcal{B}(\kappa)T} W \|_Y \leq C e^{-\nu^2 (\kappa_0(1-\delta))^2 T}\|W\|_Y.
\]
For large times $T=1$, this does not actually produce decay: $e^{-\nu^2 (\kappa_0(1-\delta))^2 T} = e^{-\nu^2 (\kappa_0(1-\delta))^2} \sim 1$. Therefore, we really do need the stronger result given in Proposition \ref{lem:intermediate-wn} to conclude that small wavenumbers $|\kappa| \leq \kappa_0$ really do give the leading order behavior of solutions.
\end{Remark}

\begin{Proof}
Let $\delta \in (0,\frac{1}{4})$ and fix $\kappa \in [\kappa_0(1-\delta), \kappa_1/\nu]$, with any fixed $\kappa_1 > \nu \kappa_0(1-\delta)$. We will study the decay of solutions to 
\[
\frac{d}{dT}\hat U = \mathcal{B}(\kappa) \hat U,
\]
with $\hat U$ and $\mathcal{B}(\kappa)$ defined in \eqref{E:defUhat} and \eqref{E:defB} using Villani's theory of hypocoercivity \cite{Villani09}. Writing this equation in components and writing $\hat u_0^s = u$ and $\hat u_m^s = v_m$ with $m = 1, 2, \dots$ for notational convenience, we have
\begin{eqnarray*}
\partial_T u &=& -\nu^2\kappa^2 u + A\rmi \kappa \sum_{m=1}^{\infty} \chi_m v_{m} \\
\partial_T v_m &=& -(\nu^2\kappa^2 + \mu_m) v_m + A\rmi \kappa \chi_m u + A\rmi \kappa \sum_{j =1}^{\infty} \chi_{m,j} v_{j}.
\end{eqnarray*}
Motivated by \cite{Villani09}, we consider the functional
\begin{equation*}
\Phi[(u, v)](T) = \zeta_0 u\bar u +  \sum_{m}\zeta_m v_m\bar v_m + 2 \mathrm{Re}\left( \rmi u \sum_{m} \sigma_m \bar v_m \right)
\end{equation*}
with $\zeta_0$, $\zeta_m$, and $\sigma_m$ to be defined below. We will show that $\dot\Phi \leq -{\color{black}\tilde M} \Phi$ for some constant {\color{black}$\tilde M$} that is independent of $\nu$ and $\kappa$, as long as $ \kappa_1/\nu \geq |\kappa| \geq \kappa_0(1-\delta)$. We will also chose $\zeta_0$, $\zeta_m$, and $\sigma_m$ so that there exist constants $c_{1,2}$ independent of $\nu$ and $\kappa$ so that {\color{black}$c_1 \|(u,v)\|_Y^2 \leq \Phi(u,v) \leq c_2 \|(u,v)\|_Y^2$}. This will imply that {\color{black}$\|(u,v)(\tau)\|_Y \leq \sqrt{c_2/c_1}e^{-\frac{1}{2}\tilde M T}$}. Undoing the scalings will then imply the decay claimed in the Proposition.

We compute
\begin{eqnarray*}
\dot\Phi &=& -2\zeta_0\nu^2\kappa^2 |u|^2 -2A\kappa\zeta_0 \mbox{Re}\left( \rmi u \sum_m  \chi_m \bar v_{m} \right) - 2 \sum_m \zeta_m (\nu^2\kappa^2 + \mu_m) |v_m|^2 + 2 A \kappa \mbox{Re} \left( \rmi u \sum_m \zeta_m \chi_m \bar v_m\right) \nonumber \\
&& \quad + 2 A\kappa \mbox{Re} \left(\rmi \sum_m \zeta_m \bar v_m \sum_{j=1}^{\infty}\chi_{m,j} v_j\right) -2\nu^2\kappa^2 \mbox{Re}\left( \rmi u \sum_m \sigma_m \bar v_m \right) - 2A\kappa \mbox{Re} \left(\sum_j \chi_j v_j \sum_m \sigma_m \bar v_m \right)  \nonumber  \\
&& \quad - 2 \mbox{Re}\left( \rmi u \sum_m \sigma_m (\nu^2\kappa^2 + \mu_m) \bar{v}_m \right) + 2 A\kappa \mbox{Re} \left( |u|^2 \sum_m \sigma_m \chi_m \right)  \\
&& \quad + 2 A\kappa \mbox{Re} \left( u \sum_m \sigma_m \sum_{j = 1}^{\infty}  \chi_{m,j} \bar v_j \right). \nonumber 
\end{eqnarray*}
Next, define
\[
\sigma_m = - \frac{c}{2A\kappa \mu_m} \chi_m, \qquad \zeta_m = \zeta_0 \quad \forall m,
\]
where $c$ is a constant to be determined. Note that this choice of $\zeta_m$ implies
\[
-2A\kappa\zeta_0 \mbox{Re}\left( \rmi u \sum_m  \chi_m \bar v_{m} \right) + 2 A\kappa \mbox{Re} \left( \rmi u \sum_m \zeta_m \chi_m \bar v_m\right) = 0
\]
Also, 
\[
2 A\kappa \mbox{Re} \left(\rmi \sum_m \zeta_m \bar v_m \sum_{j = 1}^{\infty}\chi_{m,j} v_j\right) = 0,
\]
which results from the fact that the $\chi_{m,j}$ are real and $\chi_{m,j} = \chi_{j,m}$. This follows from the fact that the eigenfunctions $\psi_j$ of the Laplacian on the cross section $\Omega$ can be chosen to be real. Therefore, we have
\begin{eqnarray*}
\dot\Phi &=& -2\zeta_0\nu^2\kappa^2 |u|^2 - 2 \zeta_0 \sum_m  (\nu^2\kappa^2 + \mu_m) |v_m|^2 +\frac{c\nu^2\kappa}{A} \mbox{Re}\left( \rmi u \sum_m \frac{1}{\mu_m} \chi_m \bar v_m \right) + c \mbox{Re} \left(\sum_j \chi_j v_j \sum_m \frac{1}{\mu_m}\chi_m \bar v_m \right)  \nonumber  \\
&& \quad + \frac{c}{A} \mbox{Re}\left( \rmi u \sum_m \frac{\chi_m}{\kappa\mu_m} (\nu^2\kappa^2 + \mu_m) \bar{v}_m \right) -c  |u|^2 |\chi|_{\mu}^2  - c \mbox{Re} \left( u \sum_m \frac{\chi_m}{\mu_m} \sum_{j = 1}^{\infty} \chi_{m,j} \bar v_{j} \right) \\
&\leq& \left[ -2\zeta_0\nu^2\kappa^2 + \frac{c\nu^2|\kappa|}{A Q_1^2}  + \frac{c}{2 A |\kappa|Q_2^2} -c |\chi|_{\mu}^2 + \frac{c}{2Q_3^2}\right] |u|^2 \\
&& \qquad  \left[ -2\zeta_0(\mu_1 + \nu^2\kappa^2) + \frac{c \nu^2 |\kappa| Q_1^2|\chi|_{\mu}^2 }{A}  + c|\chi||\chi|_\mu  + \frac{cQ_2^2}{2A |\kappa|} |\chi|^2 + \frac{cQ_3^2}{2}|\chi|_\mu^2|\chi|_{L^\infty}^2  \right]|v|^2. \\
& = & \left(I_u + II_u\right) |u|^2 + \left(I_v + II_v\right) |v|^2,
\end{eqnarray*}
where we denote $|v| = \|v\|_{\ell^2}$, {\color{black} $Q_{1,2,3}$ are constants that will be chosen later, and} where
\[
 I_u = -c |\chi|_{\mu}^2 + \frac{c}{2 A |\kappa|Q_2^2}+ \frac{c}{2Q_3^2}, \qquad 
 II_u = -2\zeta_0\nu^2\kappa^2 + \frac{c\nu^2|\kappa|}{A Q_1^2},
 \]
 and
 \[
 I_v = -2\zeta_0\mu_1 + c |\chi||\chi|_\mu + \frac{cQ_3^2}{2}|\chi|_\mu^2|\chi|_{L^\infty}^2 +  \frac{cQ_2^2}{2A |\kappa|} |\chi|^2, \qquad II_v = -2\zeta_0\nu^2\kappa^2 + \frac{c \nu^2 |\kappa| Q_1^2|\chi|_{\mu}^2 }{A}.
 \]
Recall that  $0 < \delta < 1/4$ and $|\kappa| > \kappa_0(1-\delta)$. Furthermore, let 
\[
c < \text{min } \left\{\frac{1-\delta}{|\chi|_{\mu}^2}, \frac{\mu_1}{|\chi||\chi|_\mu + |\chi|_{L^\infty}^2 + \frac{2|\chi|^2|}{3A^2 \kappa_0^2 |\chi|_\mu^2}}, A^2 \kappa_0^2 (1-\delta), \frac{12 \mu_1}{| \chi |_{\mu}^2} \right\}.
\]
We choose $\zeta_0 = 1$, $Q_1^2 = Q_2^2 = \frac{1}{A\kappa_0 |\chi|_{\mu}^2}$ and $Q_3^2 = \frac{2}{|\chi|_{\mu}^2}$. 
Then 
\begin{eqnarray*}
I_u = c |\chi|_{\mu}^2 \left( -\frac{3}{4} + \frac{\kappa_0}{2 |\kappa|} \right) \leq -\frac{c |\chi|_{\mu}^2}{12} \leq - \mu_1
\end{eqnarray*}
since $|\kappa| > \kappa_0(1-\delta)$, $0 < \delta < 1/4$, and $c < 12\mu_1/|\chi|_\mu^2$. Next, notice that the above choices imply that 
\begin{eqnarray*}
II_u = \nu^2 \left( - 2 \kappa^2 + c |\chi|_{\mu}^2 \kappa_0 |\kappa| \right) \leq - \nu^2\kappa^2,
\end{eqnarray*}
where we have used the facts that $|\kappa| > \kappa_0(1 - \delta)$ and $c < \frac{(1-\delta)}{ |\chi|_{\mu}^2}$. Similarly, 
\begin{eqnarray*}
I_v = -2\mu_1 + c \left[ |\chi||\chi|_\mu + |\chi|_{L^\infty}^2 + \frac{2|\chi|^2}{3A^2 \kappa_0^2 |\chi|_\mu^2}\right] \leq - \mu_1.
\end{eqnarray*}
Finally,  
\begin{eqnarray*}
II_v = \nu^2 \left( - 2 \kappa^2 + \frac{c }{ A^2 \kappa_0} |\kappa| \right) \leq - \nu^2 \kappa^2
\end{eqnarray*}
because $c < A^2 \kappa_0^2 (1-\delta)$. Therefore 
\begin{eqnarray*}
\dot{\Phi} \leq -(\mu_1+ \nu^2 \kappa^2) (|u|^2 + |v|^2).
\end{eqnarray*}
Also, we have that 
\begin{eqnarray*}
\Phi & \leq & \left( 1 + \frac{c}{2A|\kappa|}\right) |u|^2 + \left( 1 + \frac{c|\chi|_{\mu}^2}{2A |\kappa|}\right)|v|^2 \\
& \leq & {\color{black}\check{M}}  (|u|^2 + |v|^2),
\end{eqnarray*}
where ${\color{black}\check{M}} = 1 + \frac{A \kappa_0 }{2}\max\{1,| \chi |_{\mu}^2 \}$. As a result,
\begin{eqnarray}
\dot{\Phi} \leq - {\color{black} \tilde M} \Phi,
\end{eqnarray}
where {\color{black}$\tilde M = \mu_1/\check M$}. If we now additionally require that
\[
c \leq \mbox{min}\left\{ A\kappa_0 (1-\delta), \frac{A\kappa_0 (1-\delta)}{|\chi|_{\mu}^2}\right\}, 
\]
we find
\[
\Phi \geq \left( 1 - \frac{c}{2A\kappa_0 (1-\delta)}\right)|u|^2 + \left( 1 - \frac{c|\chi|_{\mu}^2}{2A \kappa_0 (1-\delta)}\right) |v|^2 \geq \frac{1}{2}(|u|^2 + |v|^2).
\]
Therefore, 
\[
|u(T)|^2 + |v(T)|^2 \leq 2 \Phi(T) \leq 2e^{-{\color{black} \tilde M}T}\Phi(0) \leq 4e^{-{\color{black} \tilde M}T}[|u(0)|^2 + |v(0)|^2],
\]
which completes the proof of the Proposition.
\end{Proof}


\subsection{Splitting of the semigroup} \label{sec:splitting}
The goal of this subsection is to establish the decay rates on the semigroup by splitting it as
\begin{equation}\label{eq:semSplit}
e^{B(\kappa)(T-s)} = \mathcal{E}_\mathrm{high}(\kappa, T-s) +  \mathcal{E}_\mathrm{low}(\kappa, T-s) + \mathcal{T}_N(\kappa, T-s) + \mathcal{R}_N(\kappa, T-s),
\end{equation}
where the components are defined as follows. Both $\mathcal{E}_{high, low}$ will be exponentially decaying pieces that correspond to high and low wavenumbers, respectively. The terms $\mathcal{T}_N$ and $\mathcal{R}_N$ will both correspond to the leading order eigenvalue $\lambda_0(\kappa) = -\nu_{td} \kappa^2 + \Lambda_0(\kappa)$ of $\mathcal{B}(\kappa)$, defined in Proposition \ref{prop:spec-decomp}, with $\mathcal{T}_N$ arising from the Taylor diffusion term $-\nu_{td} \kappa^2$ and $\mathcal{R}_N$ arising from the remainder $\Lambda_0(\kappa)$. 

To precisely define each term in \eqref{eq:semSplit}, first let $\psi(\kappa)$ be a smooth bump function that equals $1$ for $|\kappa| \leq \kappa_0$ and $0$ for $|\kappa| \geq 2 \kappa_0$, where {\color{black} $\kappa_0$ is a fixed small constant that is independent of $\nu$ and whose value will be specified below.} Furthermore, let $P_0(\kappa)$ be the ($\nu$-independent) projection onto the eigenspace for the leading eigenvalue $\lambda_0(\kappa)$ of $\mathcal{B}(\kappa)$, defined in Proposition \ref{prop:spec-decomp}, and let $Q_0(\kappa) = I - P_0(\kappa)$ be its complement. We can then define
\begin{eqnarray}
\mathcal{E}_\mathrm{high}(\kappa, T-s) &=& (1-\psi(\kappa)) e^{\mathcal{B}(\kappa)(T-s)} \label{def:exphigh} \\
\mathcal{E}_\mathrm{low}(\kappa, T-s) &=& \psi(\kappa)Q_0(\kappa) e^{\mathcal{B}(\kappa)(T-s)}. \label{def:explow}
\end{eqnarray}
We use a Taylor expansion to define the remaining two terms $\mathcal{T}_N$, acting on a function $\hat G(\kappa, s)$, and $\mathcal{R}_N$ as
\begin{eqnarray}
\mathcal{T}_N(\kappa, T-s)\hat G(\kappa,s) &=& e^{-\nu_{td} \kappa^2(T-s)} \sum_{\ell=0}^{N} \frac{1}{\ell!} \partial_\kappa^{\ell} \left(\psi(\kappa) P_0(\kappa) e^{\Lambda_0(\kappa)(T-s)} \hat{G}(\kappa,s) \right)|_{\kappa=0} \kappa^{\ell} \label{def:TN}\\
\mathcal{R}_N(\kappa, T-s) &=& e^{-\nu_{td} \kappa^2(T-s)} \psi(\kappa) P_0(\kappa) e^{\Lambda_0(\kappa)(T-s)}  - T_N(\kappa,T-s). \label{def:Rem}
\end{eqnarray}
With this definition, we have 
\[
\mathcal{T}_N(\kappa, T-s) + \mathcal{R}_N(\kappa, T-s) =  \psi(\kappa)P_0(\kappa) e^{\mathcal{B}(\kappa)(T-s)} = \psi(\kappa) P_0(\kappa) e^{\lambda_0(\kappa)(T-s)} =  e^{-\nu_{td} \kappa^2(T-s)} \psi(\kappa) P_0(\kappa) e^{\Lambda_0(\kappa)(T-s)}.
\]
We now obtain decay estimates on each piece of \eqref{eq:semSplit}. 


\subsubsection{Bounds on $\mathcal{E}_\mathrm{low}$} \label{sec:explow}
Before providing bounds on $\mathcal{E}_\mathrm{low}$, we first state the following lemma.
\begin{Lemma} \label{lem:L2gauss} 
Recall $\nu_{td} = \nu^2 + D_{td}$, where $D_{td} = A^2 \| \chi \|_{ \mu}^2$. Let $d > 0$ and $T > 0$. Then
\[
\| \kappa^d e^{-\nu_{td} \kappa^2 (1+T)} \|_{L^2(\mathbb{R})} \leq C (1+T)^{-\frac{d}{2}-\frac{1}{4}},
\]
where the constant $C=C(d)$ is independent of $\nu.$
\end{Lemma}
\begin{Proof}
This follows from a direct calculation; the $\nu$-independence of the constant $C$ follows from the fact that 
\[
\nu_{td}^{-1} = (\nu^2 + D_{td})^{-1} \leq D_{td}^{-1}
\]
\end{Proof}
We now prove the following lemma, which provides estimates on $\mathcal{E}_\mathrm{low}$. Recall that the norm $\| \cdot \|$ is defined in \eqref{E:Uhatnorm}.
\begin{Lemma} \label{lem:Elow} 
\begin{enumerate}
\item $\| \mathcal{E}_\mathrm{low}(\cdot,T)\hat{V}(\cdot) \| \leq C e^{-\frac{\mu_1}{2}T} \|\hat{V}(\cdot) \|$
\item $\|\int_0^T \mathcal{E}_\mathrm{low}(\cdot,T-s)\hat{F}(\cdot,s) ds \| \leq C (1+T)^{-\frac{N}{6} -\frac{1}{12}}.$
\end{enumerate}
\end{Lemma}

\begin{Proof}
By \eqref{E:Uhatnorm} and Corollary \ref{cor:explow}, we have
\begin{eqnarray*}
\|\mathcal{E}_\mathrm{low}(\cdot,T)\hat{V}(\cdot) \|^2 &=& \int_\mathbb{R} \|\mathcal{E}_\mathrm{low}(\kappa,T)\hat{V}(\kappa) \|^2_Y\rmd \kappa \\
&\leq&  \int_\mathbb{R}C e^{- \mu_1T} \|\hat{V}(\kappa) \|^2_Y\rmd \kappa = C e^{-\mu_1T} \|\hat{V}(\cdot) \|^2.
\end{eqnarray*}
This proves (i). To prove item (ii), note 
\begin{eqnarray*}
 \|\int_0^T \mathcal{E}_\mathrm{low}(\cdot,T-s)\hat{F}(\cdot,s) ds \| &=& \| \|\int_0^T \mathcal{E}_\mathrm{low}(\cdot,T-s)\hat{F}(\cdot,s) ds \|_Y \|_{L^2(\mathbb{R})} \\
 && \quad \leq \int_0^T C  e^{- \frac{\mu_1}{2}(T-s)} \| \|\hat{F}(\cdot,s)\|_Y \|_{L^2(\mathbb{R})} \rmd s.
\end{eqnarray*}
Now from Lemma \ref{lem:Fbound}, we know that 
\begin{eqnarray*}
\|\hat{F}(\kappa,s)\|_Y &\leq&  C(1+s)^{\frac{N-1}{2} - \frac{1}{2}(j+n)} e^{-\nu_{td}\kappa^2(1+s)} [|\kappa|^{N+1} + |\kappa|^{N+2} (1+s)^{1/2} + |\kappa|^{N+2} + |\kappa|^{N+3} (1+s)^{1/2}].
\end{eqnarray*}
Next, using Lemma \ref{lem:L2gauss}, we have that
\begin{eqnarray*}
 \| \|\hat{F}(\cdot,s)\|_Y \|_{L^2(\mathbb{R})} & \leq & C (1+s)^{\frac{N-1}{2} - \frac{N+1}{2} - \frac{1}{4}} (1+s)^{-\frac{1}{2}(j+n)} \\
& = & C (1+s)^{-\frac{5}{4}} (1+s)^{-\frac{1}{2}(j+n)},
\end{eqnarray*}
where the constant $C$ is independent of $\nu$.
Therefore
\begin{eqnarray*}
\|\int_0^T \mathcal{E}_\mathrm{low}(\cdot,T-s)\hat{F}(\cdot,s) ds \| &\leq& \int_0^T C  e^{- \frac{\mu_1}{2}(T-s)} (1+s)^{-\frac{5}{4}} (1+s)^{-\frac{1}{2}(j+n)} \rmd s \\
 &\leq& C (1+T)^{-\frac{1}{2}(j+n) - \frac{1}{4}} \\
 &\leq& C (1+T)^{-\frac{N}{6}-\frac{1}{12}},
\end{eqnarray*}
where the exponent in the last line follows from the fact that $N- 1 = 3j + n$, and $n \in \{0, 1, 2\}$.
\end{Proof}


\subsubsection{Bounds on $\mathcal{E}_\mathrm{high}$}
We prove the following Lemmas.
\begin{Lemma} \label{lem:Ehigh} There exist constants $C$ and $M_1$, independent of $\nu$, such that
\begin{enumerate}
\item $\| \mathcal{E}_\mathrm{high}(\cdot,T)\hat{V}(\cdot) \| \leq C e^{- M_1 T} \|\hat{V}(\cdot) \|$
\item 
\[ 
\|\int_0^T \mathcal{E}_\mathrm{high}(\cdot,T-s)\hat{F}(\cdot,s) ds \| \leq C e^{-\frac{1}{4}M_1T}.
\]
\end{enumerate}
\end{Lemma}

\begin{Proof}
We can use Proposition \ref{lem:intermediate-wn} for $\kappa_0  \leq |\kappa| \leq \frac{\kappa_1}{\nu}$ and Corollary \ref{cor:exphigh} for $|\kappa| \geq \frac{\kappa_1}{\nu}$ to find
\[
\|e^{\mathcal{B}(\kappa)T}W\|_Y \leq C e^{- M T} \|W\|_Y
\]
for all $|\kappa| \geq \kappa_0$. Therefore, we have
\begin{eqnarray*}
\|\mathcal{E}_\mathrm{high}(\cdot,T)\hat{V}(\cdot) \|^2 &=& \int_\mathbb{R} \|\mathcal{E}_\mathrm{high}(\kappa,T)\hat{V}(\kappa) \|^2_Y\rmd \kappa \\
&\leq&  \int_\mathbb{R}C e^{-2 M T} \|\hat{V}(\kappa) \|^2_Y\rmd \kappa = C e^{-2 M T} \|\hat{V}(\cdot) \|^2,
\end{eqnarray*}
which proves (i). Item (ii) follows additionally from Lemma \ref{lem:Fbound} and the estimate
\begin{eqnarray*}
&& \|\int_0^T \mathcal{E}_\mathrm{high}(\cdot,T-s)\hat{F}(\cdot,s) ds \| \leq \int_0^T C  e^{- M(T-s)} \| (1-\psi(\cdot))\|\hat{F}(\cdot ,s)\|_Y \|_{L^2(\mathbb{R})} \rmd s \\
&& \quad \leq  \int_0^T C  e^{- M(T-s)} (1+s)^{\frac{N-1-n}{3}}   \|(1-\psi(\kappa))|\kappa|^{(N+1)}e^{-\nu_{td}\kappa^2(1+s)}[1+(|\kappa| + |\kappa|^2 )(1+s)^{1/2}] \|_{L^2(\mathbb{R})} \rmd s \\
&& \quad \leq  \int_0^T C  e^{- M(T-s)} \sup_{|\kappa| \geq 2 \kappa_0} (e^{-\frac{\nu_{td}}{2} \kappa^2(1+s)})(1+s)^{\frac{N-1-n}{3}}   \||\kappa|^{(N+1)}e^{-\frac{\nu_{td}}{2}\kappa^2(1+s)}[1+(|\kappa| + |\kappa|^2 )(1+s)^{1/2}] \|_{L^2(\mathbb{R})} \rmd s \\
&& \quad \leq \int_0^T C  e^{- M(T-s)} \sup_{|\kappa| \geq 2 \kappa_0} (e^{-\frac{\nu_{td}}{2} \kappa^2(1+s)}) (1+s)^{-\frac{N}{6}-\frac{1}{2}}\rmd s \\
&& \quad = \int_0^T C  e^{- M(T-s)} e^{-2\nu_{td} \kappa_0^2(1+s)} (1+s)^{-\frac{N}{6}-\frac{1}{2}}\rmd s \\
&& \quad = \int_0^T C  e^{- M(T-s)} e^{-\nu_{td} \kappa_0^2(1+s)}e^{-\nu_{td} \kappa_0^2(1+s)} (1+s)^{-\frac{N}{6}-\frac{1}{2}}\rmd s \\
&& \quad \leq C  e^{-M T} e^{(M - \nu_{td} \kappa_0^2)T} \leq C e^{-D_{td} \kappa_0^2 T}.
\end{eqnarray*}
\end{Proof}


\subsubsection{Bounds on $\mathcal{R}_\mathrm{N}$}

In this section we prove the following lemma.
\begin{Lemma} \label{lem:RemIC}
Recall that $\mathcal{R}_\mathrm{N}(\kappa,T)$ is defined in \eqref{def:Rem}, $\| \cdot \|$ in \eqref{E:Uhatnorm}, and $\hat F$ in \eqref{E:defF}. Then
\begin{enumerate} 
\item $\| \mathcal{R}_\mathrm{N}(\cdot, T)\hat V(\cdot)\|^2 \leq C(\psi, P_0, \hat V)  T^{-\frac{N}{6}-\frac{5}{12}}$
\item $\| \int_0^T \mathcal{R}_\mathrm{N}(\cdot ,T-s) \hat{F}(\cdot ,s)ds \| \leq C  (1+T)^{-\frac{N}{6} -\frac{1}{12}}$
\end{enumerate}
for all $T > 0$, where the constant $C(\psi, P_0, \hat V)$ depends on the first $N+1$ derivatives of $\psi, P_0$, and $\hat V$. In particular, we need to require that $\|\partial_\kappa^\ell \hat V\|$ is bounded for all $0 \leq \ell \leq N+1$.
\end{Lemma}

\begin{Remark}
We can ensure that the initial condition $\hat U(\kappa, 0)$ in \eqref{E:remduhamel} has $\|\partial_\kappa^\ell \hat U(\cdot, 0)\|$ bounded for all $0 \leq \ell \leq N+1$ by requiring that the initial condition $u(x,y,z,0)$ to \eqref{E:orig-model} lies in $L^2((N+1))$.
\end{Remark}

\begin{Proof}
To estimate $\mathcal{R}_\mathrm{N}$, notice that a smooth function minus the first $N$ terms of its Taylor series can be written
\[
f(\kappa) - \sum_{j=0}^N \frac{1}{j!}f^{(j)}(0)\kappa^j = \int_0^{\kappa} \int_0^{\kappa_N} \dots \int_0^{\kappa_1} \partial_{y}^{N+1} f(y) \rmd y \rmd \kappa_1 \dots \rmd \kappa_N.
\]
Therefore, we can write
\[
\mathcal{R}_\mathrm{N}(\kappa, T-s)\hat G(\kappa, T-s) = e^{-\nu_{td} \kappa^2 T} \int_0^{\kappa} \int_0^{\kappa_N} \dots \int_0^{\kappa_1} \partial_{y}^{N+1} \left[ \psi(y)e^{\Lambda_0(y)(T-s)}P_0(y)\hat G(y, T-s) \right]\rmd y \rmd \kappa_1 \dots \rmd \kappa_N.
\]
Furthermore, the computation of the expansion of $\lambda_0(\kappa)$ that follows equation \eqref{eq:evalexp} implies that
\[
|\lambda_0(\kappa) + \nu_{td} \kappa^2| = |\Lambda_0(\kappa)| \leq C|\kappa|^3, \qquad \mbox{ for } |\kappa| \leq 2 \kappa_0
\]
for some constant $C$ that is independent of $\nu$. The $y$-derivatives in the above integral expression could fall on any of the terms in the brackets. Thus, we need to bound terms of the form
\[
\| (\partial_\kappa^{m_1}\psi(\kappa))( \partial_{\kappa}^{m_2}e^{\Lambda_0(\kappa)T})(\partial_\kappa^{m_3}P_0(\kappa))(\partial_{\kappa}^{m_4}\hat V(\kappa))\|, \qquad m_1 + m_2 + m_3 + m_4 = N+1. 
\]
Using the form of $\Lambda_0(\kappa)$, we have
\[
\partial_{\kappa}^{m_2} e^{\Lambda_0(\kappa)T} \sim (\kappa^2 T )^{\rho_1}(\kappa T)^{\rho_2} (T )^{\rho_3}(T)^{\rho_4}\cdots (T)^{\rho_{m_2}}e^{\Lambda_0(\kappa)T},
\]
where $\rho_1 + 2 \rho_2 + \dots + m_2 \rho_{m_2} = m_2$ and $\rho_i \in \{ 0, 1, \dots, m_2\}$ for all $i$. Thus, the $\rho_i$ term corresponds to $i$ derivatives falling on $\Lambda_0(\kappa)$. Therefore,
\begin{eqnarray*}
&& \| \mathcal{R}_\mathrm{N}(\cdot, T)\hat V(\cdot)\|^2 \leq \int_\mathbb{R} \left \| e^{-\nu_{td} \kappa^2 T} \int_0^{\kappa} \int_0^{\kappa_N} \dots \int_0^{\kappa_1} \partial_{y}^{N+1} \left[ \psi(y)e^{\Lambda_0(y)T}P_0(y)\hat V(y) \right] \rmd y \rmd \kappa_1 \dots \rmd \kappa_N \right \|_Y^2 \rmd \kappa \\
&& \leq \sum_{m_1+m_2+m_3+m_4 = N+1} \sup_{|\kappa| \leq 2 \kappa_0}\left\|(\partial_\kappa^{m_1}\psi(\kappa))(\partial_\kappa^{m_3}P_0(\kappa))(\partial_\kappa^{m_4}\hat V(\kappa))\right\|^2_Y \\
&& \qquad \qquad \qquad \times \int_{|\kappa| \leq 2\kappa_0}  \left( e^{-\nu_{td} \kappa^2 T}\int_0^{\kappa} \int_0^{\kappa_N} \dots \int_0^{\kappa_1} \left|\partial_y^{m_2} e^{\Lambda_0(y)T} \right| \rmd y \rmd \kappa_1 \dots \rmd \kappa_N\right)^2 \rmd \kappa \\
&&  = C(\psi, P_0, \hat V) \int_{|\kappa| \leq 2\kappa_0}  e^{-2\nu_{td} \kappa^2 T}  \\
&& \qquad  \times \left(\int_0^{\kappa} \int_0^{\kappa_N} \dots \int_0^{\kappa_1}\left| (y^2 T )^{\rho_1}(y T)^{\rho_2} (T )^{\rho_3}(T )^{\rho_4}\cdots (T )^{\rho_{m_2}}e^{\Lambda_0(y)T}  \right|\rmd y \rmd \kappa_1 \dots \rmd \kappa_N\right)^2 \rmd \kappa \\
&&  \leq C(\psi, P_0, \hat V)  T^{2(\rho_1 + \dots + \rho_{m_2})}\int_{|\kappa| \leq 2\kappa_0} e^{-2\nu_{td} \kappa^2 T}e^{2C|\kappa|^3T}|\kappa|^{2(2\rho_1+\rho_2+N+1)} \rmd \kappa. 
\end{eqnarray*}
The constant $C(\psi, P_0, \hat V)$ is determined by $\sup_{|\kappa| \leq 2 \kappa_0}\left\|(\partial_\kappa^{m_1}\psi(\kappa))(\partial_\kappa^{m_3}P_0(\kappa))(\partial_\kappa^{m_4}\hat V(\kappa))\right\|^2_Y$. The function $\psi$ and the projection $P_0$ are smooth, bounded, and independent of $\nu$, so we need not worry about derivatives that fall on them. Notice that, for $z = \kappa \sqrt{T}$, we have
\begin{eqnarray*}
\int_{|\kappa| \leq 2\kappa_0}  e^{-2\nu_{td} \kappa^2 T}e^{2C|\kappa|^3T} |\kappa|^{\rho} \rmd \kappa &=& C T^{-\frac{(\rho + 1)}{2}}\int_{|z| \leq 2\kappa_0 \sqrt{T}} |z|^{\rho}e^{-2\nu_{td} z^2} e^{2CT|z/\sqrt{T}|^3} \rmd z \\
&\leq& C T^{-\frac{(\rho + 1)}{2}}\int_{|z| \leq 2\kappa_0 \sqrt{T}} |z|^{\rho}e^{-\nu_{td} z^2}e^{-z^2\left( \nu_{td} - \frac{2C|z|}{ \sqrt{T}}\right) } \rmd z \\
&\leq& C  T^{-\frac{(\rho+1)}{2}},
\end{eqnarray*}
Note that we have used the fact that  $|z| \leq 2\kappa_0 \sqrt{T}$. Therefore, after possibly making $\kappa_0$ smaller if necessary, $\nu_{td} - \frac{2C|z|}{ \sqrt{T}} \geq 0$ . As a result
\[
\| \mathcal{R}_\mathrm{N}(\cdot, T)\hat V(\cdot)\| \leq C(\psi, P_0, \hat V) T^{\rho_1 + \dots + \rho_{m_2}} T^{-\frac{1}{4} - \frac{1}{2}(2\rho_1+\rho_2+N+1)}.
\]
Notice that
\begin{eqnarray*}
\rho_1 + \dots + \rho_{m_2} -\frac{1}{4} - \frac{1}{2}(2\rho_1+\rho_2+N+1) &=& \frac{1}{2}\rho_2 + \rho_3 + \rho_4 + \dots \rho_{m_2} -\frac{3}{4} - \frac{N}{2} \\
&\leq& \frac{1}{3}(\rho_1 + 2 \rho_2 + \dots + m_2 \rho_{m_2}) -\frac{3}{4} - \frac{N}{2} \\
&\leq& -\frac{N}{6} -\frac{5}{12}.
\end{eqnarray*}
In addition, 
\[
2\rho_1 +2 \rho_2 + 2 \rho_3 + 3 \rho_4 + \dots +(m_2-1)\rho_{m_2} \leq 2m_1 \leq 2(N+1).
\]
Therefore, we obtain
\[
\| \mathcal{R}_\mathrm{N}(\cdot, T)\hat V(\cdot)\| \leq C(\psi, P_0, \hat V)  T^{-\frac{N}{6}-\frac{5}{12}},
\]
which proves (i). To prove (ii), Lemma \ref{lem:Fbound} implies
\[
\|\hat F(\kappa, s)\|_Y \leq C  |\kappa|^{N+1} (1+s)^{ \frac{N-1}{3}}e^{-\nu_{td}\kappa^2(1+s)}[1+(|\kappa| + |\kappa|^2 )(1+s)^{1/2}].
\]
Similarly, 
\[
\|\partial_k^{m_4}\hat F(k, s)\|_Y \leq C  |\kappa|^{N+1- r_1} (1+s)^{ \frac{N-1}{3}} (\partial_\kappa^{r_2} e^{-\nu_{td} \kappa^2(1+s)})\partial_\kappa^{r_3}[1+(|\kappa| + |\kappa|^2 )(1+s)^{1/2}],
\]
where $r_1 + r_2 + r_3 = m_4$. Moreover, 
\[
(\partial_\kappa^{r_2} e^{-\nu_{td} \kappa^2(1+s)}) = \sum_{q_1 + 2q_2 = r_2} C (-\nu_{td} \kappa(1+s))^{q_1}(-\nu_{td}(1+s))^{q_2}e^{-\nu_{td} \kappa^2(1+s)} .
\]
As a result, 
\begin{eqnarray*}
&& \| \mathcal{R}_\mathrm{N}(\cdot, T-s)\hat F(\cdot, s)\|^2 \leq \\
&& \qquad \int_\mathbb{R} \left \| e^{-\nu_{td} \kappa^2(T-s)} \int_0^{\kappa} \int_0^{\kappa_N} \dots \int_0^{\kappa_1} \partial_{y}^{N+1} \left[ \psi(y)e^{\Lambda_0(y)(T-s)}P_0(y)\hat F(y, s) \right] \rmd y \rmd \kappa_1 \dots \rmd \kappa_N \right \|_Y^2 \rmd \kappa \\
&& \quad \leq \sum_{m_1+m_2+m_3+m_4 = N+1} \sup_{|\kappa| \leq 2 \kappa_0}\left\|(\partial_\kappa^{m_1}\psi(\kappa))(\partial_\kappa^{m_3}P_0(\kappa))\right\|^2_{L(Y)} \\
&& \qquad \qquad \qquad \times \int_{|\kappa| \leq 2\kappa_0} \left( e^{-\nu_{td} \kappa^2(T-s)}\int_0^{\kappa} \int_0^{\kappa_N} \dots \int_0^{\kappa_1} \left\||\partial_y^{m_2} e^{\Lambda_0(y)(T-s)} \partial_y^{m_4} \hat F(y, s)\right\|_Y \rmd y \rmd \kappa_1 \dots \rmd \kappa_N\right)^2 \rmd \kappa \\
&& \quad \leq C(\psi, P_0) \nu_{td}^{2(q_1+q_2)} (T-s)^{2(\rho_1 + \dots + \rho_{m_2})}(1+s)^{\frac{2(N-1)}{3} + 2(q_1+q_2)} \\
&& \qquad \qquad \qquad \times \int_{|\kappa| \leq 2\kappa_0} e^{-2\nu_{td} \kappa^2(T-s)}e^{-2\nu_{td} \kappa^2(1+s)}e^{2C|\kappa|^3(T-s)}|\kappa|^{2(2\rho_1+\rho_2+N+1)} |\kappa|^{2(N+1-r_1 + q_1)} \rmd \kappa \\
&& \quad \leq   C(\psi, P_0) \nu_{td}^{2(q_1+q_2)} (T-s)^{2(\rho_1 + \dots + \rho_{m_2})}(1+s)^{\frac{2(N-1)}{3} + 2(q_1+q_2)} \\
&& \qquad \qquad \times \mathrm{min}\left\{(T-s)^{-(2N+2 + 2 \rho_1 + \rho_2 - r_1 + q_1) - \frac{1}{2}}, (1+s)^{-(2N+2 + 2 \rho_1 + \rho_2 - r_1 + q_1) - \frac{1}{2}} \right\}.
\end{eqnarray*}
As a result,
\begin{eqnarray*}
&& \int_0^T \| \mathcal{R}_\mathrm{N}(\cdot, T-s)\hat F(\cdot, s)\| \rmd s \leq \\
&& \quad C(\psi, P_0)  \nu_{td}^{(q_1+q_2)}  \int_0^{T/2} (T-s)^{(\rho_1 + \dots + \rho_{m_2})}(1+s)^{\frac{(N-1)}{3} + (q_1+q_2)}  (T-s)^{-(N+1 +  \rho_1 + \rho_2/2 - r_1/2 + q_1/2) - \frac{1}{4}} \rmd s  \\
&& \quad + C(\psi, P_0)  \nu_{td}^{(q_1+q_2)}  \int_{T/2}^T (T-s)^{(\rho_1 + \dots + \rho_{m_2})}(1+s)^{\frac{(N-1)}{3} + (q_1+q_2)}  (1+s)^{-(N+1 +  \rho_1 + \rho_2/2 - r_1/2 + q_1/2) - \frac{1}{4}} \rmd s \\
&& \quad \leq C(\psi, P_0) \nu_{td}^{(q_1+q_2)} (1+T)^{-\frac{N}{6} - \frac{1}{12}} \leq C(\psi, P_0, A, \chi) (1+T)^{-\frac{N}{6} - \frac{1}{12}}.
\end{eqnarray*}
Note that we have used the fact that
\begin{eqnarray*}
&& (\rho_1 + \dots + \rho_{m_2}) + \frac{(N-1)}{3} + (q_1+q_2) - (N+1 +  \rho_1 + \rho_2/2 - r_1/2 + q_1/2 \\
&& \qquad \qquad = \left( \frac{\rho_2}{2} + \rho_3 + \dots \rho_{m_2}\right)  - \frac{2N}{3} - \frac{7}{12} + \frac{r_1}{2} + \frac{q_1}{2} + q_2 \\
&& \qquad \qquad \leq \frac{1}{3} (\rho_1 + 2 \rho_2 + \dots + m_2 \rho_{m_2})  - \frac{2N}{3} - \frac{7}{12} + \frac{1}{2}(r_1 + r_2) \\
&& \qquad \qquad \leq \frac{1}{3} m_2 + \frac{1}{2} m_4 - \frac{2N}{3} - \frac{7}{12}  \\
&& \qquad \qquad \leq \frac{1}{2}(N+1) - \frac{2N}{3} - \frac{7}{12} = -\frac{N}{6} - \frac{1}{12}.
\end{eqnarray*}
\end{Proof}


\subsubsection{Bounds on $\mathcal{T}_\mathrm{N}$}

In this section we show that the Taylor polynomial terms are actually zero. 
\begin{Lemma} \label{lem:TNzero}
Let $\mathcal{T}_\mathrm{N}(\kappa,T-s))\hat{G}(\kappa,s)$ be defined as in \eqref{def:TN}. If $s=0$ and $\hat{G}(\kappa,0)$ is an initial condition for \eqref{E:remduhamel}, or if $\hat{G}(\kappa,s) = \hat{F}(\kappa, s)$, where $\hat F$  is defined in \eqref{E:defF}, then $\mathcal{T}_\mathrm{N}(\kappa,T-s)\hat{G}(\kappa,s) = 0$ for all $\kappa$.
\end{Lemma}

\begin{Proof}
Recall that
\bea 
\mathcal{T}_\mathrm{N}(\kappa,T-s)\hat{G}(\kappa,s) = \sum_{\ell=0}^{N} \frac{1}{\ell!} \partial_\kappa^{\ell} \left( \psi(\kappa) e^{\Lambda_0(\kappa)(T-s)} P_0(\kappa) \hat{G}(\kappa,s) \right) |_{\kappa=0} \kappa^{\ell}. 
\eea 
In this expression, some derivatives fall on $\hat{G}(\kappa,s)$, but the order of these derivatives does not exceed $N$. First consider the case where $\hat{G}(\kappa,0)$ is an initial condition for \eqref{E:remduhamel}. This implies that
\[
\hat{G}(\kappa, 0) = \begin{pmatrix} \hat u_0^s(\kappa, 0) \\ \{ \hat u_n^s(\kappa, 0)\}_{n=1}^\infty \end{pmatrix}.
\]
The functions $\hat{u}_n^s$, for $n = 0, 1, \dots$ are defined via the projections in \eqref{E:deflowhigh}, the similarity variables in 
\eqref{E:cov}, and the Fourier transform. Equation \eqref{E:deflowhigh} defines $w_0^s$ and $v_n^s$ as the projections off of the first $N+1$ eigenfunctions of the operator $\mathcal{L}_{td}$. The projections onto those eigenfunctions are defined in terms of the Hermite polynomials, which implies that
\[
\int_\mathbb{R} \xi^j w_0^s(\xi, \tau) \rmd \xi = \int_\mathbb{R} \xi^j v_n^s(\xi, \tau) \rmd \xi = 0
\]
for all $\tau \geq 0$ and $j = 0, \dots, N$. Since when $\tau = 0$ we have $\xi = X$, we therefore find
\[
\partial_\kappa^j \hat u_0^s(\kappa, 0)|_{\kappa=0} = \int_\mathbb{R} \mathcal{F}^{-1}[\partial_\kappa^j \hat u_0^s(\cdot, 0)](X) \rmd X = C \int_\mathbb{R} X^j u_0^s(X,0) \rmd X = C \int_\mathbb{R} X^j w_0^s(X,0) \rmd X = 0, \quad j = 0, \dots N,
\]
where we have used $\mathcal{F}^{-1}$ to denote the inverse Fourier transform and $C$ is some constant that can be explicitly determined. Similarly, 
\[
\partial_\kappa^j \hat u_n^s(\kappa, 0)|_{\kappa=0} = C \int_\mathbb{R} X^j u_n^s(X,0) \rmd X = C \int_\mathbb{R} X^j \partial_X v_n^s(X,0) \rmd X = -j C \int_\mathbb{R} X^{j-1} v_n^s(X,0) \rmd X = 0, \quad j = 1, \dots N.
\]
When $j = 0$, the result holds because $\int \partial_\xi v_n^s(\xi,0) \rmd \xi = \int V_n(\xi, 0) \rmd \xi = 0$, where $V_n$ is defined in \eqref{E:defvnT}.

Next, consider the case where $\hat G = \hat F$. Note that $\hat{F}(\kappa,s) = \kappa^{N+1} \hat{H}(\kappa,s)$, where $\hat{H}(\kappa,s)$ is a smooth, bounded function in $\kappa$ and $s$. This fact can be seen from equation \eqref{E:defF2}. Therefore $\partial_\kappa^{\ell} \hat{G}(\kappa,s)|_{\kappa=0} = 0$ for $0 \leq \ell \leq N$.
\end{Proof}


\subsection{Proof of Proposition \ref{prop:errordecay}, and hence Theorem \ref{thm:main1}(ii)}
Recall that the goal of this chapter is to prove Proposition \ref{prop:errordecay}, which by Remark \ref{rem:thm1ii} implies Theorem \ref{thm:main1}(ii). Hence, we want to establish the estimate
\beas 
\| \hat{U}(\cdot ,T)\| \leq C(1+T)^{-\frac{N}{6} - \frac{1}{12}}.
\eeas 
Recall from \eqref{E:remduhamel} that 
\beas 
\hat{U}(\kappa,T) = e^{\mathcal{B}(\kappa)T} \hat{U}(\kappa,0) + \int_0^T e^{\mathcal{B}(\kappa)(T-s)} \hat{F}(\kappa,s) ds.
\eeas  
Using the splitting of the semigroup in \eqref{eq:semSplit} and Lemmas \ref{lem:Elow}, \ref{lem:Ehigh}, \ref{lem:RemIC}, \ref{lem:TNzero}, we have 
\begin{eqnarray*}
\|\hat U(T)\| &\leq& C \left[ e^{-\frac{ \mu_1}{2}T} + e^{-M T} + (1+T)^{-\frac{N}{6} - \frac{5}{12}}\right] \|\hat U(0)\| \\
&& \qquad + C \left[(1+T)^{-\frac{N}{6}-\frac{1}{12}} + e^{-\frac{1}{4}M T} + (1+T)^{-\frac{N}{6} - \frac{1}{12}} \right],
\end{eqnarray*}
which proves the result. 


\begin{Acknowledgment}
The authors thank Tasso Kaper and Edgar Knobloch for
useful discussions of Taylor dispersion and the anonymous referee for very helpful comments regarding the original version of this manuscript. CEW also thanks Tony Roberts
for helpful correspondence about this problem.  The work of MB is
supported in part by National Science Foundation grant DMS-1411460 and that of OC and CEW is
supported in part by National Science Foundation grant DMS-1311553.
\end{Acknowledgment}

\bibliography{taylor-dispersion-refs}
 
 
\end{document}